\documentclass[11pt,reqno]{amsart}
\usepackage{mathtools,amsthm,amssymb,amsfonts}
\usepackage{graphicx}
\usepackage{hyperref}
\usepackage{verbatim}
\usepackage{csquotes}
\usepackage{geometry}
\usepackage{color}
\usepackage{tikz-cd}
\usetikzlibrary{cd}
\usepackage{setspace}
\usepackage{todonotes}

\usepackage[shortlabels]{enumitem}
\usepackage[capitalise]{cleveref}

\usepackage{booktabs}
\usepackage{todonotes}
\usepackage{quiver}
\usepackage{comment}

\usepackage{charter,eulervm}
\usepackage[mathscr]{euscript}

\theoremstyle{plain}
\newtheorem{theorem}{Theorem}[section]
\newtheorem*{theorem*}{The Contraction Principle}
\newtheorem{corollary}[theorem]{Corollary}
\newtheorem{proposition}[theorem]{Proposition}

\newtheorem{lemma}[theorem]{Lemma}

\newtheorem{fact}[theorem]{Fact}

\theoremstyle{definition}
\newtheorem*{definition*}{Definition}
\newtheorem{definition}[theorem]{Definition}

\newtheorem{example}[theorem]{Example}

\theoremstyle{remark}
\newtheorem{remark}[theorem]{Remark}
\newtheorem*{notation*}{Notation}

\newtheorem{convention}[theorem]{Convention}

\newcommand{\map}[3]{ #1 \colon #2 \to #3 }
\newcommand{\emap}[3]{ #1 \colon #2 \hookrightarrow #3 }

\def\s{\sigma}

\def\H{\mathcal{H}}
\def\I{\mathcal{I}}

\def\B{\mathfrak{B}}
\def\D{\mathcal{D}}

\def\s{\mathfrak{s}}

\newcommand{\Frm}{\mathbf{Frm}}
\newcommand{\SLat}{\mathbf{SLat}}
\newcommand{\DSLat}{\mathbf{DSLat}}
\newcommand{\DSLatsup}{\mathbf{DSLat_{sup}}}
\newcommand{\DLat}{\mathbf{DLat}}

\newcommand{\Coz}{\mathord{\sf{Coz\,}}}

\newcommand{\ignore}[1]{}

\newcommand{\cb}{\prec\!\!\prec}
\newcommand{\rb}{\prec}

\newcommand{\meet}{\wedge}
\newcommand{\join}{\vee}

\title[Semilattice base hierarchy for frames]{Semilattice base hierarchy for frames\\and its topological ramifications}

\author{G.~Bezhanishvili}
\address{New Mexico State University, Las Cruces, NM, USA.}
\email{guram@nmsu.edu}

\author{F.~Dashiell Jr}
\address{CECAT, Chapman University, Orange, CA, USA.}
\email{dashiell@math.ucla.edu }

\author{A.~Razafindrakoto}
\address{University of the Western Cape, Cape Town, South Africa}
\email{arazafindrakoto@uwc.ac.za}

\author{J.~Walters-Wayland}
\address{CECAT, Chapman University, Orange, CA, USA.}
\email{joanne@waylands.com}

\date{}

\begin{document}



\begin{abstract}
    We develop a hierarchy of semilattice bases (S-bases) for frames. For a given (unbounded) meet-semilattice $A$, we analyze the interval in the coframe of sublocales of the frame of downsets of $A$ formed by all frames with the S-base $A$. We study various degrees of completeness of $A$, which generalize the concepts of extremally disconnected and basically disconnected frames. We introduce the concepts of D-bases and L-bases, as well as their bounded counterparts, and show how our results specialize and sharpen in these cases. Classic examples that are covered by our approach include zero-dimensional, completely regular, and coherent frames, allowing us to provide a new perspective on these well-studied classes of frames, as well as their spatial counterparts. 
\end{abstract}
\maketitle

\tableofcontents

\section{Introduction} \label{sec: intro}

For any frame $L$, a subset $A\subseteq L$ is a {\em base} of $L$ if each element of $L$ is a join of elements from $A$. As its own entity, $A$ can have structure, sometimes inherited from $L$, sometimes independent, and often a combination of both. We are interested in situations when $A$ is at least a meet-semilattice (not necessarily bounded) and the embedding $A \hookrightarrow L$ is a meet-semilattice morphism. 
Developing the theory of semilattice bases for frames was initiated by Johnstone \cite{Johnstone:1982} using sites and further developed by Ball et al. \cite{ball_extending_2014, BPWW:2016}. 
Our aim is to continue this line of research. We mainly use the language of sublocales and their corresponding nuclei instead of sites.
We pay special
attention to how the meet-semilattice sits within the frame, and to the degree of completeness of the join operation in $A$ (namely when $A$ itself is a frame or a $\sigma$-frame).
Layering additional structure onto $A$ and the embedding $A \hookrightarrow L$ gives us an interesting hierarchy of these bases. We use the following naming conventions: base $A$ is a(n) 

\begin{itemize}
    \item $S$-base: the embedding $A \hookrightarrow L$ is in $\SLat$, the category of meet-semilattices and meet-semilattice morphisms;
    \item $DS$-base: $A$ is an S-base and a distributive meet-semilattice (the embedding $A \hookrightarrow L$ is automatically in $\DSLat$, the full subcategory of distributive meet semi-lattices); 
    \item $DS_\text{sup}$-base: the embedding $A \hookrightarrow L$ is in $\DSLatsup$, the category of distributive meet-semilattices and sup-morphisms (meet-semilattice morphisms preserving existing finite joins in $A$);
    \item $L$-base: the embedding $A \hookrightarrow L$ is in $\DLat$, the category of distributive  lattices and lattice morphisms (equivalently, $A$ is a sublattice of $L$).
    \end{itemize}
Clearly each of these is a strengthening of the prior one.
 Analogous to the use of C and C$^*$ in \cite{GJ:1960} we will use $S^*$-base if the S-base $A$
 is bounded, and similarly for the other bases. 
Frames with $L^*$-bases abound and the following three examples are of particular interest.
\begin{itemize}
    \item Zero-dimensional frames: the Boolean algebra $CL$ of complemented elements.
    \item Completely regular frames: the $\sigma$-frame $\Coz L$ of cozero elements.
    \item Coherent frames: the bounded sublattice $KL$ of compact elements.
\end{itemize}

In \cref{sec: zero-dimensional,sec: cozero,sec: compact}, we concentrate on precisely these bases and use the machinery developed in the prior sections (which applies to less stringent bases) to analyze these specific categories. This yields various new characterizations of extremally disconnected, basically disconnected, and coole frames, and brings to light several interesting problems about these characterizations, depending on whether we are in the zero-dimensional or strongly zero-dimensional setting. In \cref{sec: cozero}, we introduce the notion of a cozero frame, which allows us 
to focus attention on some 
``old" topological and set-theoretic questions, as well as some interesting new open problems akin to those in \cref{sec: zero-dimensional}. 
The resulting hierarchy of various topological conditions that come into play in these contexts is summarized
in \cref{tab:ZD vs str. ZD} towards the end of \cref{sec: zero-dimensional}, and \cref{fig: Coz summary} towards the end of \cref{sec: cozero}. In \cref{sec: compact}, we apply our machinery to coherent frames \cite[Sec.~II.3]{Johnstone:1982} providing additional insight into this category.

There are interesting examples of S-bases that are not L-bases.
For example, if $L$ is semi-regular, the Booleanization $\B L$ is a DS-base 
but not generally a $DS_\text{sup}$-base.
Other examples include Dube's P-elements in completely regular frames \cite{Dube:2019} and the base associated with Martinez and Zenk's d-elements in arithmetic frames (aka~M-frames) \cite{MartinezZenk2003}. 
The distinction between $DS_\text{sup}$-bases and L-bases is more subtle.
As we will see (\cref{thm: DA in L}), if $A$ is a $DS_\text{sup}$-base of $L$, then the sublattice of $L$ generated by $A$ is exactly the distributive envelope of $A$, which then is an L-base of $L$. Nevertheless, there exist $DS_\text{sup}$-bases that are not L-bases. A rather simple example is given in \cref{ex: A neq DA}. 
Since $DS_\text{sup}$-bases play an important role in our considerations, we simplify notation and call them D-bases (see \cref{def: D-base})

For any meet-semilattice $A$, the free frame generated by $A$ is isomorphic to the frame $\D A$ of downsets of $A$ (see, e.g., \cite[Sec.~IV.2]{PicadoPultr:2011}). We give an alternate proof of \cite[Thm.~3.7]{ball_extending_2014}
that, up to isomorphism, all frames with S-base $A$ constitute an interval 
in the coframe of sublocales of $\D A$. We work with nuclei rather than coverages and the result is slightly more general in that we do not require $A$ to be bounded. 

Our results shed some light on the MacNeille completion of a meet-semilattice. 
Various degrees of distributivity of the MacNeille completion were investigated by Ern\'e \cite{Erne:1993}. In \cref{thm: MA frame} we give several equivalent conditions showing when  the MacNeille completion of a meet-semilattice is a frame akin to those studied in \cite{Cornish:1978, BPWW:2016}. \cref{prop: regular lattice,rem: sigma-regular} then give an alternate view of some of the results in \cite[Sec.~7]{ball_lindelof_2017}.

 As mentioned above, the base $A$ may have additional structure (independent of the frame -- we suggest the reader keep the idea of the Booleanization in mind) which leads us to investigate the effect of various degrees of completeness of $A$. One may ask when $A$ is a complete lattice, a complete distributive lattice (and their $\sigma$ counterparts), however in this paper we concentrate on 
when
 $A$ is itself a frame or $\sigma$-frame. This is done in \cref{sec: completeness} utilizing $D\!_\infty A$, the smallest sublocale of $\D A$ with S-base $A$, and the associated $\sigma$-frame $D\!_\sigma A$, studied in \cref{sec: S-bases}. For an S-base $A$, we obtain various characterizations of the notions of $A$-extremal, $A$-basic, and $A$-coole frames, which generalize the extremally disconnected, basically disconnected, and coole frames. 

 In \cref{sec: D-bases} we establishe that, up to isomorphism, all frames with D-base $A$ constitute the interval $[D\!_\infty A,\I D\!A]$ in the coframe of sublocales of $\D A$, where $D\!A$ is the distributive envelope of $A$ and  $\I D\!A$ is the ideal frame of $D\!A$.\footnote{We caution the reader to be aware of the difference between $D\!A$, the distributive envelope of $A$, and $\D A$, the downsets of $A$.} This provides the D-base version of \cref{thm: interval}. 
 There is then a straightforward way to obtain analogous weaker notions from \cref{sec: completeness} for D-bases, namely $D\!A$-extremal, $D\!A$-basic, and $D\!A$-coole. 

 In particular, we focus on the extreme cases at the various levels, and these often yield interesting characterizations. For example, if $A$ is bounded, we study when $D\!_\infty A$ is the Booleanization of $L$ for various L$^*$-bases $A$ of $L$. This, in particular, provides an alternate view of extremally disconnected frames (see \cref{thm: ED}) and  extremally disconnected P-frames (see \cref{prop: cozero frames}).
 
The perspective of using semilattice bases to study frames, together with our machinery, promises to be fruitful, and we conclude the paper with some future challenges and an outline of what the authors would like to accomplish.

\section{Preliminaries} \label{sec: prelims}

  Basic definitions and results that are now thought of as folklore will not be given, however we will give some definitions to establish notation and specific facts that may be less familiar but are useful to describe our results. Our primary sources are \cite{PicadoPultr:2011} for frames, \cite{Eng:1989} for topology, and \cite{AHS:2006} for category theory. We endeavor to supply accessible references for definitions and results throughout the paper, and often cite those that use similar terminology rather than the original source.

Let $\SLat$ be the category of meet-semilattices and meet-semilattice morphisms. 
We recall that a meet-semilattice $A$ is \emph{distributive} if $a\wedge b\le c$ implies there are $a',b'$ with $a\le a'$, $b\le b'$, and $a'\wedge b'=c$. 

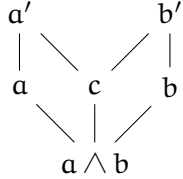
\begin{center}
    
\begin{tikzpicture}[scale=0.50]

  \node (a) at (-2,2) {$a^\prime$};

  \node (b) at (2,2) {$b^\prime$};
  \node (d) at (-2,0) {$a$};
  \node (e) at (0,0) {$c$};
  \node (f) at (2,0) {$b$};
  \node (c) at (0,-2) {$a \wedge b$};
  \draw (c) -- (d) -- (a)  -- (e) -- (b) --
  (f) -- (c) -- (e);
  \draw[preaction={draw=white, -,line width=6pt}] (a) -- (e) -- (c);
\end{tikzpicture}
\end{center}
By \cite[p.~167]{G2011}, $A$ is distributive iff the lattice of filters of $A$ is distributive (and hence a frame).
Note that Gr\"atzer works with join-semilattices and the lattices of ideals.

We point out that a weaker notion of distributivity in meet-semilattices, where the meet distributes over existing finite joins, was considered by Balbes \cite{Balbes1969} under the name of prime semilattices. Varlet \cite{Varlet1975} called these \emph{weakly distributive} semilattices, which has become standard.

Let $\Frm$ be the category of frames. The forgetful functor from $\Frm$ to $\SLat$ has a left adjoint, denoted by $\D$, that can be described by taking the collection of downsets of any meet-semilattice $A$ (see, e.g., \cite[Sec.~IV.2]{PicadoPultr:2011}). Associating with each $a\in A$ the downset ${\downarrow}a$ embeds $A$ into $\D A$, and 
any meet-semilattice morphism of $A$ to a frame $L$ factors uniquely through $\downarrow \, : A \to \D A$. 
In view of this, we will slightly abuse notation and identify $A$ with its image in $\D A$.

We will freely use the equivalence between nuclei, sublocales, and homomorphic images in $\Frm$ \cite[Ch.~III]{PicadoPultr:2011}. For a nucleus $j$, we follow \cite{MartinezZenk2003} and write $fix(j)$ to denote the sublocale of fixed points.
One of the most prominent sublocales is the Booleanization which we denote by $\B L$ with the associated nucleus  $(-)^{**}$ (which takes an element to its double pseudocomplement\footnote{We follow the standard practice of denoting the pseudocomplement of $x$ by $x^*$.}). Isbell's density theorem gives that this is the smallest dense sublocale of $L$ (see, e.g., \cite[p.~40]{PicadoPultr:2011}).

Some of this paper lives in the world of frames without the necessity of complete regularity, however this (coreflective) full subcategory, \textbf{CRFrm}, is vital to \cref{sec: zero-dimensional,sec: cozero}. 
For any frame $L$, we will denote the Boolean algebra of complemented elements of $L$ by $CL$ (thought of as the center of $L$) and the $\sigma$-frame of cozero elements by $\Coz L$. We recall (see, e.g., \cite[p.~286]{PicadoPultr:2011}) that $a \in \Coz L$ iff $a = \bigvee \{ a_n : a_n \cb a_{n+1} \}$, where $\cb$ denotes the completely below relation.

 The full subcategory of $\Frm$ consisting of compact (Lindel\"of) regular frames, \textbf{KRFrm} (\textbf{LindRFrm}), is coreflective in \textbf{CRFrm} with coreflection functor $\beta$ ($\lambda$): see, e.g., \cite[Sec.~VII.4]{PicadoPultr:2011} (\cite[Sec.~VI.7]{PicadoPultr:2021}). Recall that $\beta L$ is the frame of regular ideals of $L$ and  $\lambda L $ is the frame of $\sigma$-ideals of $\Coz L$. Thought of as sublocales, $L \subseteq \lambda L \subseteq \beta L$ (\cite[p.~476]{MaddenVermeer:1986}).
 The following facts will be useful ((1) and (2) are straightforward; for (3), see \cite{MaddenVermeer:1986} and \cite[p.~582]{BBGilmour:1996}).
\begin{fact} \label{fact: known}
For any completely regular frame $L$,
\begin{enumerate}[label=\normalfont(\arabic*), ref = \thefact(\arabic*)]
    \item \label[fact]{CL} $CL = C(\lambda L) = C(\beta L)$.
    \item \label[fact]{BL} $\B L = \B(\lambda L) =\B (\beta L)$.
    \item \label[fact]{CozL} $\Coz L = \Coz \lambda L$  but $\Coz L = \Coz \beta L$ iff $L$ is pseudocompact.  
\end{enumerate}
\end{fact}

The notions of zero-dimensionality and strong zero-dimensionality are well studied in spaces \cite{GJ:1960,Eng:1989}. Their story has also been told, albeit piecewise, in frames by various authors (see, e.g., Banaschewski and Gilmour \cite{BBGilmour:2001,BBGilmour:2001a}, Dashiell \cite{Dashiell:2011}, Dube \cite{Dube_concerning_2009}, and Marcus \cite{Nizar1998}). 
Here is a concise and simple overview told from a slightly different perspective (that fits well with our S-base philosophy).

\begin{definition}\label{defn: ZD}
$L$ is \emph{zero-dimensional} if $CL$ is a base for $L$.    
\end{definition}

Clearly every such frame is completely regular since $c \cb c$ for $c\in CL$.
Let $L_\sigma$ denote the sub-$\sigma$-frame of $L$ generated by $CL$.\footnote{Often  the center is denoted by $BL$, the ``Boolean part", and then $CL$ is used for our $L_\sigma$ (see, e.g., \cite{BB:1989,BB:2014}). Our notation helps to avoid confusing the center (when referred to as the Boolean part) with the Booleanization $\B L$ and it allows us to generalize the notation to various bases more easily.}
It is obviously a regular $\sigma$-frame, and hence contained in the largest regular sub-$\sigma$-frame, $\Coz L$. Thus, 
\[
CL \subseteq L_{\sigma} \subseteq \Coz L.
\]

\begin{definition}
    $L$ is \emph{strongly zero-dimensional} if it is zero-dimensional and every cozero element is the countable join of complemented elements, that is $\Coz L = L_{\sigma}$ (see \cite[Prop.~5.2]{Nizar1998} or \cite[Prop.~12]{Dashiell:2011}).
\end{definition}

\begin{figure}[h]
\[\begin{tikzcd}
	& L \\
	{\color{blue}{\Coz L}} && {\B L} \\
	 {\color{blue}{L_\sigma}} && {\color{blue}{(\B L)_\sigma}} \\
	& {\color{red}{CL}}
	\arrow[color={blue},hook, from=2-1, to=1-2]
	\arrow[color={green}, hook', from=2-3, to=1-2]
	\arrow[color={blue},hook, from=3-1, to=2-1]
	\arrow[color={blue},hook, from=3-3, to=2-3]
	\arrow[color={red},hook, from=4-2, to=3-1]
	\arrow[color={red},hook, from=4-2, to=3-3]
\end{tikzcd}\]
     \caption{Different bases for zero-dimensional frames}
     \label{fig: zero-dim}
 \end{figure}
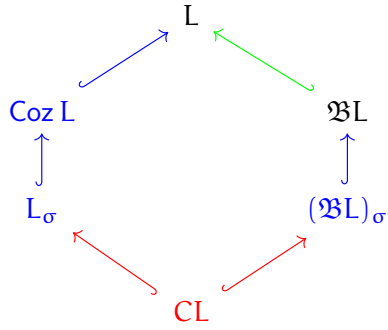
\cref{fig: zero-dim} shows the connection between different bases for zero-dimensional frames that play an important role in this paper. Black is used for objects (and morphisms) in $\Frm$, blue for  $\pmb{\sigma}\Frm$ (the category of $\sigma$-frames), red for $\DLat$, and green for $\SLat$. Note that although $\B L$ is a sublocale of $L$, it is not a subframe of $L$ since joins in $\B L$ and $L$ are different.

Let \textbf{LindZFrm} be the full subcategory of $\Frm$ consisting of Lindel\"of zero-dimensional frames. The following facts will be useful (see, e.g., \cite[Sec.~5.2]{Nizar1998} or \cite{Dashiell:2011}).

\begin{proposition}
\hfill
\begin{enumerate}[label=\normalfont(\arabic*), ref = \theproposition(\arabic*)]
    \item Every Lindel\"of zero-dimensional frame is strongly zero-dimensional.
    \item $L$ is strongly zero-dimensional iff $\lambda L$ is zero-dimensional.
    \item \label[proposition]{prop: Lind corefl3} $L$ is strongly zero-dimensional iff $\beta L$ is zero-dimensional. 
    \item \label[proposition]{prop: Lind corefl4} $\mathbf{LindZFrm}$ is a coreflective subcategory of $\mathbf{Frm}$ with coreflection map $\bigvee :\lambda_0 L\rightarrow L$ where $\lambda_0 L = \H L_\sigma$ and $\H$ is the $\sigma$-ideal functor $($see \cref{sec: S-bases}$)$.
\end{enumerate}
\end{proposition}

\begin{remark} 
\hfill
\begin{enumerate}[ref=\theremark(\arabic*)]
    \item \label[remark]{rem: strongly ZD-1} As in spaces, \cref{prop: Lind corefl3} is sometimes taken to be the definition of strongly zero-dimensional frames.
    \item The coreflection in \cref{prop: Lind corefl4} is the Lindel\"of analog of the Banaschewski compactification $\beta_0$ (where $\beta_0 L = \I CL$, the ideals of $CL$).
    \item \label[remark]{rem: strongly ZD} It is straightforward to see that $\Coz \lambda_0 L = L_\sigma$ since $\Coz \H$ is the identity on $L_\sigma$.
\end{enumerate}
\end{remark}

A natural question to ask is when $\lambda_0$ coincides with $\lambda$. The answer is rather pleasant:

\begin{corollary}
For a completely regular frame L, we have $\lambda L = \lambda_0 L$ iff $L$ is strongly zero-dimensional.
\end{corollary}
\begin{proof} 
Since $\Coz L = L_\sigma$ implies $\lambda L = \H\Coz L = \H L_\sigma = \lambda_0 L$, the reverse implication is obvious from the definition of strong zero-dimensionality. For the forward implication observe that $\Coz \lambda L = \Coz L$ (\cref{CozL}) and $\Coz \lambda_0 L = L_\sigma$ (\cref{rem: strongly ZD}). Therefore,
\[ \Coz L = \Coz \lambda L = \Coz \lambda_0 L= L_\sigma,\]
and hence $L$ is strongly zero-dimensional.
\end{proof}

\begin{remark}
This is analogous to the corresponding result for compactifications: A zero-dimensional frame $L$ is strongly zero-dimensional iff $\beta L \cong \beta_0 L$ (see, e.g., \cite[Prop.~5.2]{Nizar1998}).
\end{remark}
The methods we develop in this paper highlight the subtle difference between zero-dimensional and strongly zero-dimensional frames and emphasize several open problems related to these two concepts (see \cref{sec: zero-dimensional,sec: cozero}).

\section{S-bases} \label{sec: S-bases}

In this section we develop the basic theory of frames with a meet-semilattice base, requiring no  preservation of any additional structure of the base.
\begin{definition} \label{def: S-base}
    For any frame $L$, a subset $A \subseteq L$ is
    an \emph{S-base} if $A$ is a base and the embedding $\iota: A \hookrightarrow L$ is in $\SLat$.
     \end{definition}

With each meet-semilattice $A$, we associate two frames which will play an important role in our story. One is the frame $\D A$ of downsets of $A$ (the free frame over $A$; see \cref{sec: prelims}). The other is the frame of D-ideals of $A$, introduced by Bruns and Lakser \cite{BL70}, who showed that it is the injective hull of $A$ in $\SLat$ (see also Horn and Kimura \cite{HK71}). 
We denote it by $D\!_\infty A$. Our notation is motivated by the fact that later in the section we will have the finite and countable versions of $D\!_\infty$, which will be denoted by $DA$ and $D_\sigma A$, respectively.
We next recall some necessary definitions from \cite{BL70}.

\begin{definition} Let $A$ be a meet-semilattice. 
    \begin{enumerate}
    \item We call $G \subseteq A$ \emph{admissible} if 
    $\bigvee G$ exists and for each $a\in A$,
    $\bigvee\{ a\wedge y : y \in G\}$ also exists and 
    $
    \bigvee\{ a\wedge y : y \in G\} = a\wedge\bigvee G$.
    \item A downset $E$ of $A$ is a {\em D-ideal} if $\bigvee G\in E$ for each admissible subset $G$ of $E$. 
    \end{enumerate}
   \end{definition}
   
\begin{remark} 
\hfill
\begin{enumerate}[ref=\theremark(\arabic*)]
    \item Joins of admissible subsets of $A$ were called \emph{distributive joins} by MacNeille \cite[Def.~3.10]{MacN:1937}, the terminology that we will follow here. Ball and Pultr in \cite{ball_extending_2014} and several follow-up papers call such joins \emph{exact}.  
    We point out that MacNeille also had 
the concept of D-ideal to construct the smallest frame extension of a meet-semilattice \cite[Thm.~12.18]{MacN:1937}.\footnote{This result of MacNeille has been largely overlooked. He used the terminology \emph{lattice with completely distributive sums} for what is now called a frame.} 

\item \label[remark]{functoriality}
 In general, as is shown in \cite{AHRT2002}, taking the injective hull is not functorial. Indeed, by \cite[Thm.~2]{BL70}, an extension $e:A\hookrightarrow E$ of a meet-semilattice $A$ is essential iff $e$ preserves distributive joins and the image $e[A]$ join-generates $E$. Given that epimorphisms are onto in $\SLat$ (see \cite[Thm.~1.1]{HK71}), it is straightforward to show that essential extensions are extremal monomorphisms. Consequently,  \cite[Thm.~3.2]{AHRT2002} applies, by which $D\!_\infty$ is not a functor.
For a concrete example,
let $A$ be the Boolean algebra of finite and cofinite subsets of $\mathbb N$ and let $2$ be the two-element Boolean algebra. Clearly $D\!_\infty 2\cong 2$. On the other hand, $D\!_\infty A \cong \mathcal P(\mathbb N)$ (this, for example, can be seen by observing that $D\!_\infty A$ is isomorphic to the MacNeille completion of $A$). Let $f:A\to 2$ be the Boolean homomorphism that sends finite subsets of $\mathbb N$ to $0$ and cofinite subsets to $1$. If $g$ is an extension of $f$ to $\mathcal P(\mathbb N)$, then $g(F)=0$ for $F$ a finite subset of $\mathbb N$. On the other hand, $g(\mathbb N)=1$. Thus, $g$ cannot be a frame homomorphism. This shows that $D\!_\infty$ does not extend to a functor from $\SLat$ to $\Frm$. 
 \end{enumerate}

\end{remark}

Ball and Pultr \cite[Thm.~3.7]{ball_extending_2014} use sites to show that if $A$ is bounded then, up to isomorphism, the frames generated by $A$
form the interval $[D\!_\infty A,\D A]$ in the coframe of sublocales of 
$\D A$. 
Below we give an alternate proof 
of this result without using sites, and without the assumption that the meet-semilattice is bounded. Our proof  will yield similar intervals for more specialized bases (see, e.g., \cref{thm: D-interval}). 

We recall from \cref{sec: prelims} that we identify $A$ with its image $\{ {\downarrow}a : a \in A \}$ in $\D A$. Then $A$ 
is an S-base of $\D A$ because ${\downarrow}(a \wedge b) = {\downarrow}a \cap {\downarrow}b$ for all $a,b\in A$.
Moreover, if $A$ is an S-base for $L$, the freeness of $\D A$ gives an onto frame homomorphism $\D A \longrightarrow L$, yielding that $L$ is isomorphic to a sublocale of $\D A$.

\begin{convention} 
If $A$ is an S-base of a frame $L$, we identify $L$ with its isomorphic copy in $\D A$ and think of it as a sublocale of $\D A$.
\end{convention}

We first show that $D\!_\infty A$ is a sublocale of $\D A$ by giving an explicit description of the nucleus $\delta$ on $\D A$ such that $D\!_\infty A=fix(\delta)$. For $E \in \D A$ let
\[
\delta(E) = \left\{ \bigvee G : G \subseteq E \mbox{ is admissible}\right\}.
\]
\begin{lemma} \label{lem: D-ideal nucleus}
$\delta$ is a nucleus on $\D A$ such that $fix(\delta)=D\!_\infty A$.
\end{lemma}

\begin{proof}
    By \cite[Lem.~3]{BL70}, $\delta(E)$ is a D-ideal. This yields that $\delta$ is a closure operator on $\D A$ and that $fix(\delta)=D\!_\infty A$. It is left to show that for any downsets $E$ and $H$, $\delta(E)\cap \delta(H)\subseteq\delta(E\cap H)$. For this it is sufficient to show that $\delta(E)\cap H\subseteq\delta(E\cap H)$ (see, e.g., \cite[p.~47]{PicadoPultr:2011}). Let $a \in \delta(E)\cap H$. Then $a=\bigvee G$ for some admissible $G\subseteq E$ and $a\in H$. Therefore, $G\subseteq E\cap{\downarrow}a \subseteq E \cap H$, and so $a\in\delta(E\cap H)$. 
\end{proof}

\begin{theorem} \label{thm: interval}
The collection of frames with S-base $A$ coincides with the interval $[D\!_\infty A,\D A]$ in the coframe of sublocales of $\D A$.

\end{theorem}

\begin{proof} 

By \cref{lem: D-ideal nucleus}, $D\!_\infty A$ is a sublocale of $\D A$, and it is clear that $A$ is an S-base for $D\!_\infty A$. 
It remains to show that $D\!_\infty A$ is itself a sublocale of any frame $L$ with S-base $A$.
For any such $L$, let $P$ be the sublocale $L\cap D\!_\infty A$. Since $A$ is an S-base for both $L$ and $D\!_\infty A$, it is also an S-base for $P$. We thus have the following commutative diagram, where $f,g,p$, and $q$ are the corresponding quotients (black indicates objects and morphisms in $\Frm$ and green is used for $\SLat$):

\[\begin{tikzcd}
	{\textcolor{green}A} \\
	& {\D A} && {D\!_\infty A} \\
	\\
	& L && P
	\arrow["d"', color={green}, tail, from=1-1, to=2-2]
	\arrow["f", two heads, from=2-2, to=2-4]
	\arrow["e", color={green}, curve={height=-12pt}, tail, from=1-1, to=2-4]
	\arrow["g"', two heads, from=2-2, to=4-2]
	\arrow["q"', two heads, from=4-2, to=4-4]
	\arrow["p", two heads, from=2-4, to=4-4]
        \end{tikzcd}\]

Since $D\!_\infty A$ is an injective hull of $A$, it is an essential extension of $A$. Therefore, since $pe:A\to P$ is a one-to-one meet-semilattice homomorphism, $p:D\!_\infty A\to P$ must be one-to-one. Consequently, $D\!_\infty A=P$, and hence $D\!_\infty A\subseteq L$.
\end{proof}

\begin{remark}\hfill
\begin{enumerate}
    \item In the above proof, $P$ is the pushout of $L$ and $D\!_\infty A$ in $\Frm$.
\item If the meet-semilattice $A$ has a least element, then $D\!_\infty A$ is a dense sublocale of $L$, and thus $\B L \subseteq D\!_\infty A \subseteq L$. In certain cases (for example, see \cref{prop: cozero frames}), it is interesting to consider properties of $L$ that are equivalent to the extreme cases, that is when $D\!_\infty A$ is the Booleanization or the entire frame. This is a theme we explore throughout the paper. 
\end{enumerate}
\end{remark}

Next, for any $L$ with S-base $A$,
we give an explicit description of the nucleus corresponding to the sublocale $D\!_\infty A$ of $L$ . Define $k_{\scriptscriptstyle L}^{\scriptscriptstyle A}:L\to L$ by 
\[
k_{\scriptscriptstyle L}^{\scriptscriptstyle A}(x) = \bigwedge_{\scriptscriptstyle L} \{a\in A : x\le a\}.
\]
(When $A$ or $L$ are clear from the context, we suppress the subscript/superscript or both.) 

It is straightforward to see that $k$ is a closure operator on $L$ (for idempotency observe that $x\leq a$ iff $k_{\scriptscriptstyle L}^{\scriptscriptstyle A}(x) \leq a$). Therefore, the fixpoints $fix(k)$ form a complete lattice.
     Moreover, $A$ is both join-dense and meet-dense in $fix(k)$, and thus $fix(k)$ is isomorphic to the MacNeille completion of $A$ (see, e.g, \cite[p.~237]{BD:1974}), which we will denote by $MA$.

In general, $k$ may not be a nucleus. 
But we can use it to define a nucleus on $L$ as follows. Recall (see, e.g., \cite{Simmons1978}) that for each $c\in L$, $w_c$ defined by $w_c(x)=(x\!\to\! c)\!\to \!c$ is a nucleus on $L$, where $\to$ is the Heyting implication (relative pseudocomplement) on $L$.\footnote{Observe that $w_0(x)$ is the double pseudocomplement of $x$, that is $w_0(x)=x^{**}$.} Then, following \cite[Sec.~5]{BH14}, define
\[
j_{\scriptscriptstyle L}^{\scriptscriptstyle A}(x)=
\bigwedge_{\scriptscriptstyle L} \{w_c(x) : c \in fix(k_{\scriptscriptstyle L}^{\scriptscriptstyle A}) \}.
\]
(Again, when $A$ or $L$ are clear from the context, we suppress the subscript/superscript or both.)

\begin{lemma} \label{lem: j and k}\hfill
       \begin{enumerate}[label=\normalfont(\arabic*), ref = \arabic*]
        \item $j$ is a nucleus on $L$.
        \item \label{lem: j and k:2} $x \le j x \le k x$ for each $x\in L$, and hence $fix(k)\subseteq fix(j)$.
    \end{enumerate}
\end{lemma}

\begin{proof}
    (1) This is obvious since each $w_c$ is a nucleus and the meet of nuclei is a nucleus.

    (2) Let $x\in L$. Since $j$ is a nucleus, $x\le j x$. Moreover, $$w_{k x}(x) = (x \!\to\! k x) \!\to\! k x = 1 \to kx = k x.$$ Since $kx\in fix(k)$, by the definition of $j$ we have $j x\le k x$. Thus, $fix(k) \subseteq fix(j)$.
\end{proof}

 Since normal ideals are intersections of principal ideals (see, e.g., \cite[p.~235]{BD:1974} or \cite[p.~63]{G2011}), and $k_{\scriptscriptstyle  \D A}^{\scriptscriptstyle { A}}(E) = \bigcap \{ {\downarrow} a: E \subseteq {\downarrow} a\}$, we obtain:
\begin{lemma} \label{lem: normal ideal}
   If $A$ is a meet-semilattice, then $fix(k_{\scriptscriptstyle  \D A}^{\scriptscriptstyle { A}})$ is the collection of normal ideals of $A$.
\end{lemma}

To determine $fix (j_{\scriptscriptstyle  \D A}^{\scriptscriptstyle { A}})$, we recall that in $\D A$ the Heyting implication is calculated by 
\[
E \rightarrow H = \{ a \in A : a \meet x \in H\ \forall x \in E\}.
\]
We use this in the proof of the next proposition, which generalizes \cite[Prop.~8.12]{BH14}.

\begin{proposition} \label{prop: D-ideals}
   If $A$ is a meet-semilattice, then $\delta = j_{\scriptscriptstyle  \D A}^{\scriptscriptstyle { A}}$.
   \end{proposition}

\begin{proof}
Let $E\in\D A$ and $a\in\delta(E)$. Then $a=\bigvee G$ for some admissible $G\subseteq E$. It is sufficient to show that $a\in (E\to H)\to H$ for each $H\in fix(k_{\scriptscriptstyle  \D A}^{\scriptscriptstyle { A}})$. Let $x\in E\to H$. Then $x\wedge y\in H$ for all $y\in G$. Since $G$ is admissible, 
\[
x\wedge a = x \wedge \bigvee G = \bigvee \{x\wedge y : y\in G\} \in H
\]
because $H$ is a normal ideal by \cref{lem: normal ideal}, hence $H$ is closed under existing joins. Therefore, $a \in (E\to H)\to H$, and so $a\in j_{\scriptscriptstyle  \D A}^{\scriptscriptstyle { A}}(E)$. Thus, $\delta(E) \subseteq j_{\scriptscriptstyle  \D A}^{\scriptscriptstyle { A}}(E)$.

For the reverse inclusion, let $a\in j_{\scriptscriptstyle  \D A}^{\scriptscriptstyle { A}}(E)$. We let $E_a = \{ a\wedge e : e \in E\}$. Then $E_a\subseteq E$. We show that $a=\bigvee E_a$. Clearly $a$ is an upper bound of $E_a$. Let $x$ be another upper bound of $E_a$. Then $a\wedge e \in E_a \subseteq {\downarrow}x$ for all $e\in E$. Therefore, $a\in E\to{\downarrow}x$. Since $a\in j_{\scriptscriptstyle  \D A}^{\scriptscriptstyle { A}}(E)$, we also have $a\in (E\to{\downarrow}x)\to{\downarrow}x$, so $a\le x$. Thus, $a=\bigvee E_a$. We next show that $E_a$ is admissible. Let $c\in A$. Then $c\wedge a\in j_{\scriptscriptstyle  \D A}^{\scriptscriptstyle { A}}(E)$, so by the above,
\[
c\wedge\bigvee E_a = c\wedge a = \bigvee E_{c\wedge a} = \bigvee \{ c\wedge a\wedge e : e \in E \} = \bigvee \{ c\wedge d : d\in E_a \}.
\]
Consequently, $E_a \subseteq E$ is admissible and $a = \bigvee E_a$, yielding that $a \in \delta(E)$.
\end{proof}

As an immediate consequence of \cref{lem: D-ideal nucleus,prop: D-ideals} we obtain:

\begin{corollary}  \label{cor: fixpoints}
 If $A$ is a meet-semilattice, then $fix(j_{\scriptscriptstyle  \D A}^{\scriptscriptstyle { A}} )= D\!_\infty A$.
\end{corollary}

\cref{prop: D-ideals} provides another description of the nucleus $\delta$, which will be central to the rest of the paper. Yet another description can be given using relative annihilators. We recall that a {\em relative annihilator} of $A$ is a downset of the form 
\[
\langle a,b \rangle = \{ x\in A : a\wedge x\le b\}.
\]
In other words, $\langle a,b \rangle = {\downarrow}a \to {\downarrow}b$ in the frame $\D A$.

\begin{lemma} \label{lem: relative annihilators}
For each $E\in\D A$ we have $\delta(E)=\bigcap\{ \langle a,b \rangle : E \subseteq \langle a,b \rangle \}$.
\end{lemma}

\begin{proof}
    Let $x\in\delta(E)$ and $E\subseteq \langle a,b \rangle$. Then $x=\bigvee G$ for some admissible $G\subseteq E$. We have
   
     \[
    a \wedge x = a \wedge \bigvee G = \bigvee \{ a\wedge y: y\in G \} \le b
    \]
    since $E \subseteq \langle a,b \rangle $. Therefore, $x\in \langle a,b \rangle$, and so $\delta(E)\subseteq\bigcap\{ \langle a,b \rangle : E \subseteq \langle a,b \rangle \}$.

    For the reverse inclusion, let $x\in \langle a,b \rangle$ for each $\langle a,b \rangle \supseteq E$. Then $x\wedge a\le b$. As in the proof of \cref{prop: D-ideals}, let ${E_x = \{ x\wedge e : e \in E\}}$.
We show that $x=\bigvee E_x$. Let $u$ be an upper bound of $E_x$. Then $x\wedge e\le u$ for each $e\in E$. Therefore, $E\subseteq \langle x,u \rangle$, so $x\in \langle x,u \rangle$ by assumption, and hence $x\le u$. Arguing as in the proof of \cref{prop: D-ideals} yields that $E_x$ is admissible. Thus, $\bigvee E_x\in\delta(E)$, and so $x\in\delta(E)$.  
\end{proof}

Returning to \cref{cor: fixpoints}, we show that if $A$ is an S-base of $L$, then the fixpoints of $j^{\scriptscriptstyle { A}}_{\scriptscriptstyle {L}}$ are exactly the D-ideals of $A$. For this we need the following lemma, in the proof of which we repeatedly use that if $L$ is a sublocale of $N$, then meets of subsets of $L$ are the same in $L$ and $N$ (see, e.g., \cite[p.~26]{PicadoPultr:2011}).

\begin{lemma} \label{lem: sublocales}
    Let $A$ be an S-base of both $L$ and $N$ with $L$ a sublocale of $N$.
    \begin{enumerate}[label=\normalfont(\arabic*), ref = \thelemma(\arabic*)]
        \item \label[lemma]{lem: sublocale1} 
        $k_{\scriptscriptstyle  L}^{\scriptscriptstyle A}(x)= k_{\scriptscriptstyle  N}^{\scriptscriptstyle A}(x)$ for each $x \in L$.
        \item \label[lemma]{lem: sublocale2}
        $fix(k_{\scriptscriptstyle  L}^{\scriptscriptstyle A}) = fix(k_{\scriptscriptstyle  N}^{\scriptscriptstyle A})$.
        \item $j_{\scriptscriptstyle  L}^{\scriptscriptstyle A}(x)= j_{\scriptscriptstyle  N}^{\scriptscriptstyle A}(x)$ for each $x \in L$.
    \end{enumerate}
\end{lemma}
\begin{proof}
    (1) Let $x\in L$. Using that $A\subseteq L$ and $L$ is a sublocale of $N$, we have
    \[
   k_{\scriptscriptstyle  L}^{\scriptscriptstyle A}(x) = \bigwedge_L \{ a \in A : x \le a \} = \bigwedge_N \{ a \in A : x \le a \} = k_{\scriptscriptstyle  N}^{\scriptscriptstyle A}(x).
   \]
   
    (2) If $x \in fix(k_{\scriptscriptstyle  L}^{\scriptscriptstyle A})$, then clearly $x \in L$, and so $x=k_{\scriptscriptstyle  L}^{\scriptscriptstyle A}(x)=k_{\scriptscriptstyle  N}^{\scriptscriptstyle A}(x)$ by (1). 
   But if $x \in fix(k_{\scriptscriptstyle  N}^{\scriptscriptstyle A})$, then again $x\in L$ by the definition of $k_{\scriptscriptstyle  N}^{\scriptscriptstyle A}$ because $A\subseteq L$ and meets of subsets of $A$ are the same in $L$ and $N$. Therefore, using 
   (1) again,  $x=k_{\scriptscriptstyle  N}^{\scriptscriptstyle A}(x)=k_{\scriptscriptstyle  L}^{\scriptscriptstyle A}(x)$. Thus, $fix(k_{\scriptscriptstyle  L}^{\scriptscriptstyle A}) = fix(k_{\scriptscriptstyle  N}^{\scriptscriptstyle A})$. 

    (3) Let $x\in L$. By (2), $fix(k_{\scriptscriptstyle  N}^{\scriptscriptstyle A}) \subseteq L$. Therefore, by \cite[p.~26]{PicadoPultr:2011}, $x\to_N c \in fix(k_{\scriptscriptstyle  N}^{\scriptscriptstyle A}) \subseteq L$ for each $c \in fix(k_{\scriptscriptstyle  N}^{\scriptscriptstyle A})$. Thus, $x\to_N c = x \to_L c$. The same reasoning yields 
    \[
    (x\to_N c)\to_N c = (x\to_L c)\to_L c.
\] 
    Consequently,
    \begin{eqnarray*}
     j_{\scriptscriptstyle  N}^{\scriptscriptstyle A}(x) &=&
     \bigwedge_N \{ (x\to_N c)\to_N c : c \in fix(k_{\scriptscriptstyle  N}^{\scriptscriptstyle A}) \} \\
     &=& \bigwedge_L \{ (x\to_L c)\to_L c : c \in fix(k_{\scriptscriptstyle  L}^{\scriptscriptstyle A}) \} \\
     &=& j_{\scriptscriptstyle  L}^{\scriptscriptstyle A}(x).
    \end{eqnarray*}
\end{proof}

As an immediate consequence, we obtain:

\begin{theorem} \label{thm: char}
 For any $L$ with S-base $A$, $j_{\scriptscriptstyle  L}^{\scriptscriptstyle A}$ is the restriction of $j_{\scriptscriptstyle  \D A}^{\scriptscriptstyle A}$ to $L$ and $fix(j_{\scriptscriptstyle  L}^{\scriptscriptstyle A})= D\!_\infty A$.
   \end{theorem}
\begin{proof}
Since $L$ is (isomorphic to) a sublocale of $\D A$ (see \cref{thm: interval}), it follows from \cref{lem: sublocales} that $j_{\scriptscriptstyle  L}^{\scriptscriptstyle A}$ is the restriction of $j_{\scriptscriptstyle  \D A}^{\scriptscriptstyle A}$ to $L$ and that $fix(j_{\scriptscriptstyle  L}^{\scriptscriptstyle A}) = fix(j_{\scriptscriptstyle  \D A}^{\scriptscriptstyle A})$. However, $fix(j_{\scriptscriptstyle  \D A}^{\scriptscriptstyle { A}} )= D\!_\infty A$ by \cref{cor: fixpoints}, and hence $fix(j_{\scriptscriptstyle  L}^{\scriptscriptstyle A})= D\!_\infty A$.
\end{proof}

We conclude this section by applying our results to the MacNeille completion $MA$ of a meet-semilattice $A$. 
We view $MA$ as the complete lattice of normal ideals of $A$ (see, e.g., \cite[Sec.~XII.2]{BD:1974}). It is well known 
that in general MacNeille completions do not satisfy any distributivity laws. Several authors have investigated when $MA$ is a frame (see, e.g., \cite{Cornish:1978,Erne:1993,BPWW:2016,ball_lindelof_2017}). We utilize the results of this section to obtain the following characterizations (parts of which are proved in the above cited papers, however our approach is different).

\begin{theorem} \label{thm: MA frame}
    For a meet-semilattice $A$, the following are equivalent.
    \begin{enumerate}[label=\normalfont(\arabic*), ref = \arabic*]
        \item $MA$ is a frame.
        \item \label{thm: MA frame:2} $MA=D\!_\infty A$.
        \item \label{thm: MA frame:3} $k^A_{\D A}=j^A_{\D A}=\delta$.
        \item \label{thm: MA frame:4} $\langle a,b \rangle$ is a normal ideal for each $a,b\in A$.
    \end{enumerate}
\end{theorem}

\begin{proof}
$(1)\Rightarrow(2)$ Since $A$ is an S-base of $MA$, as sublocales $D\!_\infty A \subseteq MA \subseteq \D A$ by \cref{thm: interval}. But, every normal ideal is a D-ideal, so $MA\subseteq D\!_\infty A$, and thus $MA = D\!_\infty A$.

$(2)\Rightarrow(1)$ This is obvious.

$(2)\Leftrightarrow(3)$ This follows from \cref{lem: normal ideal,prop: D-ideals,cor: fixpoints}.

    $(2)\Rightarrow(4)$ By \cref{lem: D-ideal nucleus,lem: relative annihilators}, each relative annihilator $\langle a,b \rangle$ is in $D\!_\infty A$. Therefore, each $\langle a,b \rangle$ is a normal ideal. 

    $(4)\Rightarrow(2)$ By \cref{lem: j and k}(\ref{lem: j and k:2}), $fix(k^A_{\D A})\subseteq fix(j^A_{\D A})$. Since each $\langle a,b \rangle$ is a normal ideal and the intersection of normal ideals is normal, it follows from \cref{lem: relative annihilators,lem: normal ideal} that $fix(\delta) \subseteq fix(k^A_{\D A})$. By  \cref{prop: D-ideals}, $\delta =j^A_{\D A}$, giving the desired result.
\end{proof}

\begin{remark}
Meet-semilattices satisfying the last condition in the above theorem are called proHeyting in \cite{BPWW:2016}, where $\langle a,b \rangle$ is denoted by $b{\downarrow}a$.   
\end{remark}

We end the discussion of when the MacNeille completion is a frame by deriving some of  \cite[Sec. 7]{ball_lindelof_2017}. We recall that a distributive lattice $A$ with $1$ is \emph{subfit} (or \emph{conjunctive}) if $a\not\le b$ implies the existence of $c\in A$ such that $a\vee c=1$ and $b\vee c\ne 1$.
This notion has a long history; see \cite[p.~21]{PicadoPultr:2021}.

\begin{proposition} \label{prop: regular lattice}
    If $A$ is a subfit lattice, then $MA$ is a frame.  
\end{proposition}

\begin{proof}
By \cref{thm: MA frame}(\ref{thm: MA frame:2}), it is sufficient to show that every D-ideal is normal. Let $E$ be a D-ideal and $a$ a lower bound of the set of upper bounds of $E$. It is enough to show that $a=\bigvee E_a$ since the end of the proof of \cref{prop: D-ideals} then yields that $E_a$ is an admissible subset of $E$. Let $b$ be an upper bound of $E_a$. If $a\not\le b$, then since $A$ is a subfit lattice, there is $c \in A$ such that $a\vee c=1$ and $b\vee c\ne 1$. For $e\in E$, we have 
\[
e=e\wedge(a\vee c)=(e\wedge a)\vee(e\wedge c) \le b\vee c.
\]
Therefore, $b\vee c$ is an upper bound of $E$, so $a\le b\vee c$ by assumption. Thus, $1=a\vee c\le b\vee c$. The obtained contradiction proves that $a\le b$, and hence $a=\bigvee E_a$. 
\end{proof}

\begin{remark}\label{rem: sigma-regular}
    Recall that a bounded distributive lattice $A$ is \emph{regular} if $a\not\le b$ implies the existence of $c\in A$ such that $c \rb a$ and $c \not\le b $, and each regular lattice is subfit (see, e.g., \cite[p.~360]{ball_lindelof_2017}). Thus, by \cref{prop: regular lattice}, if $A$ is a regular lattice, $MA$ is a frame. In particular, $MA$ is a frame for each regular $\sigma$-frame.
   \end{remark}

\section{Degrees of completeness} \label{sec: completeness}

A meet-semilattice may have additional structure (possibly independent of the frame it generates).
 In this section we investigate when  the S-base $A$ is a frame or $\sigma$-frame in its own right.
 When this is the case, $A$ is automatically bounded and thus we make that assumption throughout this section:

 \begin{definition}
    We call an S-base $A$ of a frame $L$ an \emph{S$^*$-base} if $A$ is bounded.\footnote{We are using the analogy of C vs C$^*$; see \cite{GJ:1960}.}
\end{definition}

\begin{remark}
    Since $A$ join-generates $L$,  $A$ is bounded iff $A$  contains the bounds of~$L$.
\end{remark}

We first consider when $A$ is a frame. 
To give the reader some context, recall extremal disconnectivity (see \cref{def: ED}), which motivates our next definition.
 \begin{definition} \label{def: A-extremal}
    A frame $L$ with S$^*$-base $A$ is  \emph{$A$-extremal} if $A$ is a frame.
\end{definition}

We note that in $A$-extremal frames, $A$ is not usually a subframe of $L$. In fact, this holds iff $A=L$.
To characterize $A$-extremal frames, we recall that an order-preserving map $f:P\to Q$ between two posets has a left adjoint iff each set $\{ x \in P : q\le f(x) \}$ has a least element. If $P,Q$ are meet-semilattices, then we call the left adjoint {\em left exact} if it preserves binary meets. Finally, we recall from the previous section that we identify $A$ with its image in $D\!_\infty A$, and that $k,j$ are abbreviations for $k_L^A$, $j_L^A$, respectively. 

\begin{theorem} \label{thm: A extremal}
Let $A$ be an S$^*$-base of $L$. The following are equivalent.
\begin{enumerate}[label=\normalfont(\arabic*), ref = \arabic*]
    \item $L$ is $A$-extremal.
    \item $A=D\!_\infty A$.
    \item $k$ is a nucleus and $A= fix(k)$.
    \item $k=j$ and $A= fix(j)$.
    \item The embedding $A \hookrightarrow L$ has a left exact left adjoint.
\end{enumerate}
\end{theorem}

\begin{proof}
(1)$\Leftrightarrow$(2), (4)$\Rightarrow$(3), and (3)$\Rightarrow$(1) are on the surface.

(3)$\Rightarrow$(4): Since $k$ is a nucleus, we have $k=\bigwedge\{w_d: d\in fix(k)\}$ (see, e.g., \cite[Lem.~7]{Simmons1978}). But the right hand side is $j$ by definition and so, $k=j$, showing that $A= fix(j)$.

(5)$\Rightarrow$(3): If the embedding $\emap{\iota}{A}{L}$ has a left exact left adjoint $\ell$, then the composition $\iota \circ\ell$
is a nucleus on $L$ and $A$ is its fixpoints. 
Since $(\iota \circ \ell)(x)=\displaystyle\bigwedge_L\{a \in A : x \le a\}$, we obtain that $k=\iota \circ\ell$ showing that $k$ is a nucleus on $L$ with $A = fix(k)$. 

(1)$\Rightarrow$(5): Define $\map{\ell}{L}{A}$ by $\ell(x)=\displaystyle\bigvee_{A} \{ a \in A : a \le x \}$. To see that $\ell$ is left adjoint to the embedding $\iota:A\to L$, it is enough to show that $\ell(x)$ is the least element of $\{b \in A : x \le b\}$. 
If $x\not\le\ell(x)$, then since $A$ join-generates $L$, there is $a\in A$ with $a\le x$ and $a\not\le\ell(x)$. Because this contradicts the definition of $\ell(x)$, we obtain that $\ell(x)\in \{ b \in A : x \le b \}$. Moreover, for each $a,b\in A$, if $a\le x \le b$, then $a\le b$, and hence $\ell(x)\le b$. Thus, $\ell(x)$ is the least element of $\{b\in A : x \le b \}$.  

To see that $\ell$ is left exact, let $x,y\in L$. Since $A$ is a frame, we have
\begin{eqnarray*}
\ell(x) \wedge \ell(y) &=& \bigvee_{A}\{ a \in A : a \le x \} \wedge \bigvee_{A}\{ b \in A : b \le y \} \\
&=& \bigvee_{A}\{ c \in A : c \le x \wedge y \}  \\
&=& \ell(x\wedge y).
\end{eqnarray*}

\end{proof}

We now consider the situation when the S$^*$-base $A$ is a $\sigma$-frame. 

\begin{definition} \label{def: Lsigma}
    For  any bounded meet-semilattice $A$ and $L\in[D\!_\infty A,\D A]$, we let \emph{$L_{\sigma\!_A}$} be the sub-$\sigma$-frame of $L$ generated by $A$. \end{definition}
    
We simplify the notation by dropping the $A$ and writing $L_\sigma$ when the setting is unambiguous.\footnote{This is consistent with the notation in \cref{sec: prelims} where $L_\sigma$ denotes the sub-$\sigma$-frame generated by~$CL$.} When $L=D\!_\infty A$, we write $D\!_\sigma A$ for $L_\sigma$. In that case, we also let $D\!A$ be the bounded sublattice of $D\!_\infty A$ generated by $A$. We think of $D$ and $D\!_\sigma$ as the finite and countable versions of $D\!_\infty$. As we will see, $D\!A$ plays a significant role in \cref{sec: D-bases}.

\begin{proposition}\label{prop: Dinfty and Dsigma}
    For any bounded meet-semilattice $A$, 
$$D\!_\infty A= D\!_\infty D\!A \text{  and  } D\!_\sigma A= D\!_\sigma D\!A.$$
\end{proposition}

\begin{proof}
We have that $A$ is an S-base of $D\!_\infty D\!A$ and $D\!A$ is an S-base of $D\!_\infty A$. Therefore, \cref{thm: interval} implies that as sublocales, $D\!_\infty A \subseteq D\!_\infty D\!A$ and $D\!_\infty D\!A \subseteq D\!_\infty A$. Thus, $D\!_\infty A = D\!_\infty D\!A$.
That $D\!_\sigma A= D\!_\sigma D\!A$ is obvious since $D\!A \subseteq D\!_\sigma A$.
\end{proof}

For the next definition, the reader is advised to have the notion of basic disconnectivity for guidance (see \cref{def: BD}). 

\begin{definition}
    Any frame $L$ with S$^*$-base $A$ is \emph{$A$-basic} if $A$ is a $\sigma$-frame.
\end{definition}
\color{black}

\begin{theorem} \label{thm: A basic}
Let $A$ be an S$^*$-base of $L$. The following are equivalent.
\begin{enumerate}[label=\normalfont(\arabic*), ref = \arabic*]
    \item $L$ is $A$-basic.
    \item $A=D\!_\sigma A$.
   \item The embedding $A \hookrightarrow L_\sigma$ has a left exact left adjoint.
    \item There is a nucleus $\nu$ on $L_\sigma$ such that $A = fix (\nu)$.
\end{enumerate}
\end{theorem}

\begin{proof}
(1)$\Leftrightarrow$(2) and (4)$\Rightarrow$(1) are on the surface.

(3)$\Rightarrow$(4): Let $\emap{\iota}{A}{L_\sigma}$ be the embedding and $\map{\ell}{L_\sigma}{A}$ its left exact left adjoint. Then $\nu := \iota\circ\ell$ is a nucleus on $L_\sigma$ whose fixpoints is $A$.

(1)$\Rightarrow$(3): Let $x \in L_\sigma$. Since $L_\sigma$ is the sub-$\sigma$-frame of $L$ generated by $A$, there is a countable family $\{a_n\}\subseteq A$ such that $x=\displaystyle\bigvee_L a_n$. Define $\ell:L_\sigma\to A$ by $\ell(x)=\displaystyle\bigvee_A a_n$. To see that $\ell$ is well defined, let $x=\bigvee_L a_n = \bigvee_L b_n$. Then $a_n \le \bigvee_L b_n \le \bigvee_A b_n$, so $a_n\le\bigvee_A b_n$ for each $n$. Therefore, $\bigvee_A a_n \le \bigvee_A b_n$. By symmetry we get equality.  To see that $\ell$ is left adjoint to the embedding $A\hookrightarrow L_\sigma$, let $x\in L_\sigma$ and $a\in A$. Then 
\[
\ell(x)\le a \Leftrightarrow \bigvee_{A} a_n \le a \Leftrightarrow a_n\le a \ \forall n \Leftrightarrow \bigvee_L a_n \le a \Leftrightarrow x\le a.
\]
Finally, to see that $\ell$ is left exact, let $x,y\in L_\sigma$ and assume that $x=\displaystyle\bigvee_L a_n$ and $y=\displaystyle\bigvee_L b_n$. W.l.o.g., $\{ a_n\}$ and $\{ b_n\}$ may be taken to be chains, and hence $x\wedge y = \displaystyle\bigvee_L (a_n \wedge b_n)$. Since $A$ is closed under finite meets,
\begin{eqnarray*}
\ell(x) \wedge \ell(y) &=& \left(\bigvee_{A} a_n\right) \wedge \left(\bigvee_{A} b_n\right) = \bigvee_{A} (a_n \wedge b_n) 
= \ell(x\wedge y). 
\end{eqnarray*}
\end{proof}

The next obvious criterion to consider is when $D\!_\sigma A$ is itself already a frame. For this we require the following:

\begin{lemma} \label{lem: sigma frame quotients}
Let $A$ be an S$^*$-base for $L$ and $N$. Any frame quotient $\map{f}{L}{N}$ that fixes $A$ restricts to a $\sigma$-frame quotient $\map{f}{L\!_\sigma}{N\!_\sigma}$.
   \end{lemma}
\begin{proof} Since $A$ is an S$^*$-base for both $L$ and $N$, the restriction $\map{f}{L\!_\sigma}{N\!_\sigma}$ is a well-defined $\sigma$-frame map. To show it is onto,  
    let $x \in N_\sigma$. Then 
    $x=\displaystyle\bigvee_{\scriptscriptstyle N} a_n$ for some $ a_n \in A$. Now, let 
   $y=\displaystyle\bigvee_L a_n$. Then $y\in L_\sigma$ and  
    \[
    f(y)=f\left(\bigvee_L a_n\right) = \bigvee_{\scriptscriptstyle N} f(a_n) = \bigvee_{\scriptscriptstyle N} a_n = x. 
    \]
\end{proof}

As an immediate consequence we obtain:

\begin{theorem}\label{thm:fixpoints Dsigma}\label{thm: sigma restriction}
For any $L$ with S$^*$-base $A$, the restriction of ${\map{j_{\scriptscriptstyle L}}{L} { D\!_\infty A}}$ to $L_\sigma$ is a $\sigma$-frame homomorphism from $L_\sigma$ onto $D\!_\sigma A$. Thus, 
    the restriction of $j_{\scriptscriptstyle L}$ is a nucleus on $L_\sigma$ whose fixpoints are $D\!_\sigma A$.
\end{theorem}
\begin{corollary}
For a bounded meet-semilattice $A$, $D\!_\sigma A$ is a $\sigma$-frame quotient of $(\D A)\!_\sigma$.
\end{corollary}

The terminology in the next definition will be justified before \cref{defn:coole}.

\begin{definition}
    Any frame $L$ with S$^*$-base $A$ is \emph{$A$-coole} if $D\!_\sigma A=D\!_\infty A$. 
\end{definition}

\begin{theorem} \label{thm: Dsigma=Dinfty}
Let $A$ be an S$^*$-base of $L$. The following are equivalent.
\begin{enumerate}[label=\normalfont(\arabic*), ref = \arabic*]
     \item $L$ is $A$-coole. 
     \item $D\!_\sigma A$ is a frame.
     \item $fix(j)=j[L_\sigma]$.
    \item $fix(j) \subseteq L_\sigma$.
\end{enumerate}
\end{theorem}

\begin{proof}
    (1)$\Leftrightarrow$(2) and (4)$\Rightarrow$(3) are on the surface and 
(3)$\Rightarrow$(4) follows from \cref{thm: sigma restriction}.
    
    (3)$\Rightarrow$(1): It is enough to observe that $D\!_\sigma A=j[L_\sigma]$ (\cref{thm: sigma restriction}) and $D\!_\infty A = fix(j)$ (\cref{thm: char}). 
    
    (1)$\Rightarrow$(3): It is enough to show that $fix(j) \subseteq j[L_\sigma]$. Let $x\in fix(j)$. Then $x \in D\!_\infty A$ by \cref{thm: char}, so $x \in D\!_\sigma A$ by (1). Therefore, $x=\displaystyle\bigvee_{D\!_\infty A} a_n$ for some countable family $\{ a_n \} \subseteq A$. Let $y= \displaystyle\bigvee_L a_n$. Then $y \in L_\sigma$ and $jy=x$ (because $j$ is a $\sigma$-frame homomorphism by \cref{thm: sigma restriction} and $ja_n=a_n$ for each $n$). Thus, $x\in j[L_\sigma]$.
\end{proof}
\begin{remark}\label{rem:Dsigma extremal}
Clearly $L$ is $A$-coole iff $L$ is $D\!_\sigma A$-extremal.
 For a non-trivial example of this phenomenon, one needs a frame $L$ with $S^*$-base $A$ such that $D\!_\sigma A$ is neither $A$ nor $L$. Such an example will be provided in \cref{ex: Q1}.
\end{remark}

\section{D-bases and L-bases}\label{sec: D-bases}

In this section we consider two special classes of S-bases, which we term D-bases and L-bases; see the introduction. As we will see, they have various advantages over S-bases. 
Consider $\B L$:  a frame $L$ is called \emph{semi-regular} when $\B L$ is a base, in which case $\B L$ is an 
S$^*$-base. While $\B L$ is a complete Boolean algebra, the embedding need not preserve even finite joins in $\B L$. 
This happens iff $x^{**} \join y^{**} = (x^{**} \join y^{**})^{**}$ for all $x,y \in L$ (i.e., iff $L$ is extremally disconnected; see the paragraph after \cref{def: ED}). 
This motivates us to look at morphisms that preserve existing finite joins more closely.

\begin{definition} \cite{BJ2011} 
Let $A,B$ be distributive meet-semilattices. We call a meet-semilat\-tice morphism $\map hAB$ a \emph{sup-morphism} if it preserves existing finite joins in $A$.     
\end{definition}

Let $\DSLatsup$ be
 the category of distributive meet-semilattices and sup-morphisms.
 Let also $\DLat$ be the category of distributive lattices and lattice morphisms. 
 Clearly $\DLat$ is a full subcategory of $\DSLatsup$.
It follows from the result of Cornish and Hickman \cite{CH1978} that the embedding ${\DLat}\hookrightarrow{\DSLatsup}$ has a left adjoint, which associates with each distributive meet-semilattice $A$ its distributive envelope $D\!A$. 

\begin{remark}
 We emphasize the functoriality of $D$ compared with the non-functoriality of $D\!_\infty$. The example given in \cref{functoriality}
 shows that $D\!_\infty$ does not extend to a functor from ${\DSLatsup}$ to $\Frm$ either.
\end{remark}

We next describe an alternative construction of $D\!A$ which is more convenient for our purposes. This construction is taken from \cite{BJ2011} (see also \cite{HP2008} where the authors work in the signature of join-semilattices).

Call a filter $P$ of $A$ \emph{prime} if it is a meet-prime element of the lattice ${\sf Filt}(A)$ of filters of $A$ (that is, $F\cap G \subseteq P$ implies $F\subseteq P$ or $G\subseteq P$ for all $F,G\in{\sf Filt}(A)$). Let $X_A$ be the poset of prime filters of $A$ and let $\s:A\to\wp(X_A)$ be the Stone map
\[
\s(a)=\{ P \in X_A : a \in P \}.
\]  
Then $\s$ embeds $A$ into $\wp(X_A)$, and $D\!A$ is the sublattice of $\wp(X_A)$ generated by the image $\s[A]$. In other words, $S\in D\!A$ iff $S=\displaystyle\bigcup_{i=1}^n \s(a_i)$ for some $a_1,\dots,a_n\in A$. 
Note that if $A$ has $1$, then $\s(1)=X_A$ and if $A$ has $0$, then $\s(0)=\varnothing$.

\begin{lemma} \label{lem: prop of D(M)}
\hfill
\begin{enumerate}[label=\normalfont(\arabic*), ref = \thelemma (\arabic*)]
    \item \cite[Lem.~3.2]{BJ2011} Let $A$ be a distributive meet-semilattice and $a_1,\dots,a_n,b \in A$. Then 
\[
\bigcap_{i=1}^n {\uparrow}a_i \subseteq {\uparrow}b \mbox{ iff } \s(b) \subseteq \bigcup_{i=1}^n \s(a_i).
\]
\label[lemma]{lem: prop of D(M)-1}
\item \cite[p.~90]{BJ2011} The embedding $\mathfrak s : A \to D\!A$ is a sup-morphism. \label[lemma]{lem: prop of D(M)-2}
\item \cite[Prop.~3.6]{BJ2011} Let $A,B$ be distributive meet-semilattices and $\map hAB$ a sup-morphism. Then there exists a unique lattice morphism $\map {D\!h}{D\!A}{D\!B}$ such that ${D\!h\circ\s_A=\s_B\circ h}:$
\[\begin{tikzcd}
	A && D\!A \\
	\\
	B && D\!B
	\arrow["h"', from=1-1, to=3-1]
	\arrow["{\s\!_A}", hook, from=1-1, to=1-3]
	\arrow["{\s\!_B}"', hook, from=3-1, to=3-3]
	\arrow["D\!h", dashed, from=1-3, to=3-3]
\end{tikzcd}\]
\label[lemma]{lem: prop of D(M)-3}
\end{enumerate}

\end{lemma}

The next theorem, which is a consequence of \cite[Thm 1.3]{CH1978}, provides the necessary machinery to realize that under the correct conditions on an S-base $A$ of a frame $L$, $D\!A$ remains within $L$. We give an alternative proof that utilizes \cref{lem: prop of D(M)}.

\begin{theorem} \label{thm: char of D(M)}
Let $A$ be a distributive meet-semilattice. Suppose that there is an embedding $\map hAB$ of $A$ into a distributive lattice $B$ such that $h$ is a sup-morphism. Then $D\!A$ is isomorphic to the sublattice of $B$ generated by $h[A]$. 
\end{theorem}

\begin{proof}
Since $D\!B$ is isomorphic to $B$, it follows from \cref{lem: prop of D(M)-3} that there is a lattice morphism $\map{f}{D\!A}{B}$ such that $f\circ\s=h$.
\[\begin{tikzcd}
	&& {D\!A} \\
	\\
	A && B
	\arrow["\s", from=3-1, to=1-3]
	\arrow["h", from=3-1, to=3-3]
	\arrow["f", dashed, from=1-3, to=3-3]
\end{tikzcd}\]
Let $E$ be the sublattice of $B$ generated by $h[A]$. Then each element of $E$ is a finite join of elements from $h[A]$. Therefore, $f$ maps $D\!A$ onto $E$. We show that $f$ is one-to-one. For this it is sufficient to show that if $x\not\subseteq y$ in $D\!A$, then $f(x)\not\le f(y)$ in $B$. Suppose that $x=\displaystyle\bigcup_{i=1}^n \s(a_i)$ and $y=\displaystyle\bigcup_{j=1}^m \s(b_j)$ for some $a_1,\dots,a_n,b_1,\dots,b_m\in A$. Then there is $a_i$ such that $\s(a_i)\not\subseteq \textbf{}\displaystyle\bigcup_{j=1}^m \s(b_j)$. By \cref{lem: prop of D(M)-1}, $\displaystyle\bigcap_{j=1}^m {\uparrow}b_j \not\subseteq {\uparrow}a_i$. Therefore, there is $c\in A$ with $b_j\le c$ for each $b_j$ but $a_i\not\le c$.  
If $f(x)\le f(y)$, then $f\left(\displaystyle\bigcup_{i=1}^n \s(a_i)\right) \le f\left(\displaystyle\bigcup_{j=1}^m \s(b_j)\right)$, and so $\displaystyle\bigvee_{i=1}^n f(\s(a_i)) \le \displaystyle\bigvee_{j=1}^m f(\s(b_j))$ . Therefore, $h(a_i) \le \displaystyle\bigvee_{j=1}^m h(b_j) \le h(c)$ (because $b_j \le c$ implies $h(b_j)\le h(c)$), and hence $a_i \le c$ because $h$ is an embedding. The obtained contradiction proves that $f(x)\not\le f(y)$. Consequently, $f$ is one-to-one, and so $D\!A$ is isomorphic to $E$. 
\end{proof}

\begin{definition} \label{def: D-base}
    For any frame $L$, we call a base $A$ a(n)
    \begin{enumerate}
        \item {\emph{D-base}} if 
      the embedding $\iota: A \hookrightarrow L$ is in $\DSLatsup$;
       \item {\emph{L-base}} if 
      the embedding $\iota: A \hookrightarrow L$ is in $\DLat$.
       \end{enumerate}
\end{definition}

In view of this definition, as an immediate consequence of \cref{thm: char of D(M)}, we obtain:

\begin{theorem} \label{thm: DA in L}
For any D-base $A$ of $L$, $D\!A$ is isomorphic to the sublattice of $L$ generated by $A$.
\end{theorem}

\begin{convention} \label{conv: identification}
    Henceforth we will identify $A$ with its image $\s[A]$ in $D\!A$, and $D\!A$ with its isomorphic copy in $L$. So we will think of $A$ as a sub-meet-semilattice of $D\!A$ and of $D\!A$ as a sublattice of $L$. Then \cref{lem: prop of D(M)-1} takes on the following form, where the join is taken in $D\!A$:
    \[
\bigcap_{i=1}^n {\uparrow}a_i \subseteq {\uparrow}b \mbox{ iff } b \le \bigvee_{i=1}^n a_i.
\]
\end{convention}

\begin{remark} 
Observe that if $A$ is a D-base, then $D\!A$ is an L-base;
 and for any L-base $A$, obviously $A$ is a D-base and $A=DA$. However, there exist D-bases that are not L-bases (see \cref{ex: A neq DA}). Because of the close connection between the two, we concentrate on D-bases but our machinery also applies to L-bases.  
\end{remark}

Analogous to the downset functor, the forgetful functor from $\Frm$ to $\DLat$ has a left adjoint, the ideal functor $\I$. (If 
 the lattice does not have a least element, we define $\varnothing$ to be the least element of the ideal frame.)
The lattice may then be identified with the collection of principal ideals (see, e.g., \cite[Sec.~II.3]{Johnstone:1982}).
Composing $\I$ with $D$ yields: 

\begin{proposition} \label{prop: JD}
  $\map{\I D}{ \DSLatsup}{ \Frm}$ is the left adjoint to the forgetful functor 
    from $\Frm$ to $\DSLatsup$.
\end{proposition}

 \begin{remark} Recalling that a {\em Frink ideal} of a distributive meet-semilattice $A$ is a nonempty subset $I$ of $A$ such that $S^{ul} \subseteq I$ for each finite $S \subseteq I$, where $S^u$ is the upper bounds and $S^l$ the lower bounds of $S$ in $A$, we have that 
     $\I D\!A$ is nothing more than the frame of Frink ideals of $A$ \cite[Cor.~4.4]{BJ2011}. 
  \end{remark}

\begin{proposition}\label{prop:Largest D-basis}
    For a distributive meet-semilattice $A$, $A$ is a D-base for  $\I D\!A$, and any $L$ for which $A$ is a D-base is isomorphic to a sublocale of $\I D\!A$.
\end{proposition}
\begin{proof} Since ${A \hookrightarrow D\!A}$ is a sup-morphism and ${\downarrow}: D\!A \hookrightarrow \I D\!A$ is a lattice morphism, we have that
${A \hookrightarrow \I D\!A}$ is a sup-morphism. Since $A$ join-generates $D\!A$ and 
$D\!A$ join-generates $\I D\!A$, we obtain that
\mbox{$\{ {\downarrow}a \in \I D\!A: a \in A\}$}, which is isomorphic to $A$, join-generates $\I D\!A$. 
For $L$ with $A$ as D-base, $\I D\!A \rightarrow L$ exists by \cref{prop: JD}, and is onto since $A$ join-generates $L$. Thus, $L$ is isomorphic to a sublocale of $\I D\!A$.
\end{proof}

\begin{convention}
    From now on we identify $L$ with its isomorphic copy in $\I D\!A$ and view $L$ as a sublocale of $\I D\!A$.
\end{convention}
 Notice that the embedding \mbox{$A \hookrightarrow \D A $} is rarely a sup-morphism and thus $A$ is rarely a D-base for $\D A$. Indeed, if there are $a,b\in A$ such that $a\vee b$ exists in $A$ and $a,b\ne a\vee b$, then the embedding does not preserve this join. In fact we have the following.
\begin{proposition} For a distributive meet-semilattice $A$, $A \hookrightarrow \D A $ is a sup-morphism iff $\D A = \I D\! A$.
\end{proposition}
\begin{proof}
   One implication is obvious. For the other, observe that $\I D\! A$ is always a sublocale of $\D A$ by \cref{thm: interval}. Moreover, if the embedding $A \hookrightarrow \D A $ is a sup-morphism, then $A$ is a D-base for $\D A$, and applying \cref{prop:Largest D-basis} gives that $\D A$ is a sublocale of $\I D\! A$.
\end{proof}
We now have that $\I D\!A$ is the largest frame for which $A$ is a D-base. We next show that  $D\!_\infty A$ is still the least frame with $A$ as a D-base.

\begin{proposition}
    
 \label{prop: D vs Dinfty}
   If $A$ is a distributive meet-semilattice, then 
    A is a D-base for $D\!_\infty A$.
      
\end{proposition}

\begin{proof} It is enough to observe that the embedding $A \hookrightarrow D\!_\infty A$
     is a sup-morphism, which is clear since  $\mathfrak s : A \hookrightarrow D\! A$ is a sup-morphism (see \cref{lem: prop of D(M)-2}) and $D\! A \hookrightarrow D\!_\infty A$ is a lattice morphism.   
     \end{proof}   
 \begin{remark}
     Another way to prove \cref{prop: D vs Dinfty} is to observe that the embedding ${A \hookrightarrow D\!_\infty A}$ preserves all distributive joins (see, e.g., \cite[Cor.~2]{BL70}), and since $A$ is distributive, all finite joins in $A$ are distributive (see \cite[Thm.~4.1]{Balbes1969}).
 \end{remark}   

In general $D\!A$ is not a frame. For example, if $A$ is a distributive lattice then $A = D\!A$. Thus, if $A$ is a distributive lattice that is not a frame then neither is $D\!A$. 
On the other hand, \cref{prop: D vs Dinfty} does imply the following.
\begin{corollary} 
For any distributive meet-semilattice $A$, $D\!A$ is a frame iff  $D\!A = D\!_\infty A$.
\end{corollary}
In \cref{thm: DA extremal} we explore other equivalences of this phenomenon. We now establish the ``D-base" version of \cref{thm: interval}. 

\begin{theorem} \label{thm: D-interval}
The collection of frames with D-base $A$ coincides with the interval $[D\!_\infty A,\I D\!A]$ in the coframe of sublocales of $\D A$.
\end{theorem}

\begin{proof} 
The proof of \cref{thm: interval} works by replacing $\D$ with $\I D$ and using \cref{prop:Largest D-basis}. 
\end{proof}

The closure operators and nuclei defined in \cref{sec: S-bases} clearly still have meaning here.

\begin{proposition}
 For any D-base $A$ of $L$, $fix({j_{\scriptscriptstyle  L}^{\scriptscriptstyle {D\!A}}}) = D\!_\infty A$.  
\end{proposition}
\begin{proof}
It is enough to observe that $D\!A$ is an S-base of $L$ and apply \cref{thm: char,prop: Dinfty and Dsigma}.
\end{proof}

This together with \cref{thm: char} yields:

\begin{corollary}\label{cor: j nucleus} 
For any D-base $A$ of $L$, the nuclei $j_{\scriptscriptstyle  L}^{\scriptscriptstyle A}$ and $j_{\scriptscriptstyle  L}^{\scriptscriptstyle {D\!A}}$ are the same.
\end{corollary}

In fact a stronger result for the closure operators $k_L^A$ and $k_L^{D\!A}$ holds:

\begin{proposition} \label{prop: kA = kDA}
    For any D-base $A$ of $L$,
$k_L^A = k_L^{D\!A}$. 
\end{proposition}

\begin{proof}
    It suffices to show that $k_L^A x = x$ for each $x \in D\!A$. Since $A$ is a base of $L$, it is enough to show that for each $a\in A$, if $a\le k_L^A x$ then $a\le x$. Suppose $a \in A$ with  $a\le k_L^A x$. Now $x=\bigvee_i b_i$ for some $b_1,\dots,b_n\in A$, and so, by the formula in \cref{conv: identification}, to conclude that $a \leq \bigvee_i b_i$ we show that  $\bigcap{\uparrow} b_i \subseteq {\uparrow}a$ in $A$:  let $c\in \bigcap{\uparrow} b_i$,  then $b_i\le c$ for each $i$. Therefore, $x\le c$, so $k_L^A x\le c$, and hence $a\le c$. 
   \end{proof}

We point out that \cref{prop: kA = kDA} yields an alternative proof of \cref{cor: j nucleus}. Recalling from \cref{sec: S-bases} that the fixpoints of $k_L^A$ are (isomorphic to) $MA$ and the fixpoints of $k_L^{D\!A}$ are (isomorphic to) $M(D\!A)$, we also obtain:

\begin{corollary} 
The MacNeille completions of $A$ and $D\!A$ are isomorphic.    
\end{corollary}

We conclude this section by considering the consequences of \cref{sec: completeness} applied to $D\!A$ instead of $A$. For this, 
  as before, it is convenient to assume that $A$ is bounded. 
  
\begin{definition}
    A D-base is called a \emph{D$^*$-base} if $A$ is bounded; and similarly one obtains L$^*$-bases.\footnote{The $L^*$ in $L^*$-base should not be confused with $\{ a^* : a \in L \}$ (which is the Booleanization of a frame $L$).}
\end{definition}
Examples of frames with bounded L-bases include zero-dimensional, completely regular, and coherent frames (see \cref{sec: zero-dimensional,sec: cozero,sec: compact}).

In general, $D\!A$ properly contains $A$, thus the notions of $A$-extremal and $A$-basic 
weaken to those of $D\!A$-extremal and $D\!A$-basic. On the other hand, as we will see, the notions of $A$-coole and $D\!A$-coole always coincide.

\begin{definition}
    Any frame $L$ with D$^*$-base $A$ is \emph{$D\!A$-extremal} if $D\!A$ is a frame.
\end{definition}

We have a characterization of $D\!A$-extremal frames that is similar to the characterization of $A$-extremal frames given in \cref{thm: A extremal}. We skip the proof because it is essentially the same as that of the aforementioned theorem. We 
remind the reader that $D\!_\infty A =D\!_\infty D\!A$ (\cref{prop: Dinfty and Dsigma}), $j=j_L^{A}= j_L^{D\!A}$ (\cref{cor: j nucleus}), and $k=k_L^{A}= k_L^{D\!A}$ (\cref{prop: kA = kDA}).

\begin{theorem} \label{thm: DA extremal}
Let $A$ be a D$^*$-base of $L$. The following are equivalent.
\begin{enumerate}[label=\normalfont(\arabic*), ref = \arabic*]
    \item $L$ is $D\!A$-extremal.
    \item $D\!A=D\!_\infty A$.
    \item $k$ is a nucleus and $D\!A= fix (k)$.
    \item $k = j$ and $D\!A=fix(j)$.
    \item The embedding $D\!A \hookrightarrow L$ has a left exact left adjoint.
\end{enumerate}
\end{theorem}
\begin{center}
    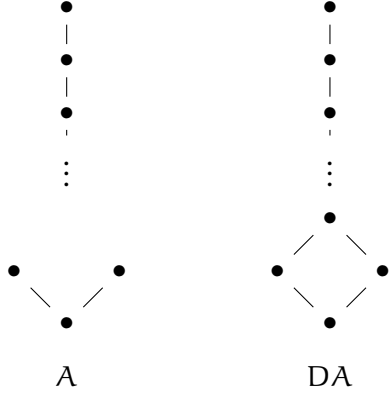
\begin{figure}[h!]
\begin{tikzpicture}[scale=.7]
  \node (max1) at (2,7) {$\bullet$};
  \node (max2) at (7,7) {$\bullet$};
  \node (a1) at (2,6) {$\bullet$};
 \node (a2) at (7,6) {$\bullet$};
  \node (b1) at (2,5) {$\bullet$};
 \node (b2) at (7,5) {$\bullet$};
 \node (c1) at (2,4) {$\vdots$};
 \node (c2) at (7,4) {$\vdots$};
  \node (d) at (7,3) {$\bullet$};
  \node (e1) at (1,2) {$\bullet$};
  \node (e2) at (6,2) {$\bullet$};
 \node (f1) at (3,2) {$\bullet$};
  \node (f2) at (8,2) {$\bullet$};
  \node (min1) at (2,1) {$\bullet$};
  \node (min2) at (7,1) {$\bullet$};
  \node  at (2,0) {$A$};
  \node  at (7,0) {$DA$};
  \draw (min1) -- (f1);
  \draw (min1) -- (e1);
   \draw (min2) -- (f2) --(d);
  \draw (min2) -- (e2) --(d);
  \draw (max1) -- (a1) -- (b1) --(c1);
  \draw (max2) -- (a2) -- (b2) --(c2);
 \end{tikzpicture}
        \caption{$A$ and $D A$}
        \label{fig: A vs DA}
    \end{figure}
\end{center}
\begin{example} \label{ex: A neq DA}
Clearly if $A$ is a frame, then so is $D\!A$ because in that case $A$ is (isomorphic to) $D\!A$. On the other hand, the diagram in \cref{fig: A vs DA} exhibits a bounded distributive meet-semilattice $A$ such that $A$ is not a lattice (see \cite[Example~3.11]{BJ2011}), but $D\!A$ is a frame. Hence, the notions of $A$-extremal and $DA$-extremal frames do not coincide. 
\end{example}

\begin{definition}
    Any frame $L$ with D$^*$-base $A$ is \emph{$D\!A$-basic} if $D\!A$ is a $\sigma$-frame.
\end{definition}
Essentially the same proof as that of \cref{thm: A basic} yields the following characterization of $D\!A$-basic frames.
We follow the convention established after \cref{def: Lsigma} and write $L_\sigma$ for $L_{\sigma_{D\!A}}$.
\begin{theorem} \label{thm: DA basic}
Let $A$ be a $D^*$-base of $L$. The following are equivalent.
\begin{enumerate}[label=\normalfont(\arabic*), ref = \arabic*]
    \item $L$ is $D\!A$-basic.
    \item $D\!A=D\!_\sigma A$.
    \item The embedding $D\!A \hookrightarrow L_\sigma$ has a left exact left adjoint.
    \item There is a nucleus $\nu$ on $L\!_\sigma$ such that $D\!A = fix(\nu)$.
\end{enumerate}
\end{theorem}

\begin{remark}
   On the other hand, a frame is $A$-coole iff it is $D\!A$-coole because $D\!_\sigma A=D\!_\sigma DA$ and $D\!_\infty A=D\!_\infty D\!A$ (see \cref{prop: Dinfty and Dsigma}). In turn, this is equivalent to being $D\!_\sigma A$-extremal (see \cref{rem:Dsigma extremal}).
In addition, when $A$ is an L$^*$-base, the distinction between $A$-extremal and $D\!A$-extremal (as well as between $A$-basic and $D\!A$-basic) vanishes.  
\end{remark}

\section{The \texorpdfstring{$L^*$}--base of complemented elements} \label{sec: zero-dimensional}

Observe that the Boolean algebra of complemented elements $CL$ is a base for $L$ iff $L$ is zero-dimensional, in which case $CL$ is clearly an $L^*$-base.
We now consider the effect of the previous sections when $L$ is zero-dimensional. First we look at the two extreme cases; when $D\!_\infty CL$ is $\B L$ or $L$ (recall that we identify $D\!_\infty CL$ with a sublocale of $L$; see \cref{thm: interval}). The next proposition shows that the one extreme case always holds.

\begin{proposition} \label{prop: L**}
For $L$ zero-dimensional, $D\!_\infty CL = \B L$.
\end{proposition}

\begin{proof}
Since $CL$ is bounded, $D\!_\infty CL$ is a dense sublocale of $L$, and hence $\B L \subseteq D\!_\infty CL$. But $CL$ is a $D$-base for $\B L$, so by \cref{thm: D-interval}, $D\!_\infty CL \subseteq \B L$.
\end{proof}

\begin{remark} \label{rem: j=k=regular}
The above clearly means that $jx=x^{**}$ for each $x\in L$. (Recall that we abreviate $j^{CL}$ to $j$ and similarly for $k$.) Moreover, since 
$x^*=\bigvee\{a\in CL : a\le x^* \}$, we have
\begin{eqnarray*}
 x^{**}=\bigwedge\{a^* : a\in CL, \ x^{**}\le a^*\} = \bigwedge\{c : c\in CL, \ x\le c\}=kx.  
\end{eqnarray*}
Thus, $fix(k)=fix(j)=\B L$. 
This can also be seen by observing that since $CL$ is a Boolean algebra, its MacNeille and D-ideal completions are isomorphic. But the former is isomorphic to the fixpoints of $k$ (see \cref{lem: normal ideal}) and the latter to the fixpoints of $j$ (see 
\cref{cor: fixpoints}).
\end{remark}

For the other extreme case we consider when $D\!_\infty CL=L$. 

\begin{proposition}
    For $L$ zero-dimensional, $D\!_\infty CL= L$ iff $L$ is a complete Boolean algebra $($aka a Boolean frame$)$.
\end{proposition}

\begin{proof}
    By \cref{prop: L**}, $D\!_\infty CL= L$ iff $\B L=L$, which means that $L$ is a Boolean frame.
\end{proof}
We now consider what \cref{thm: A extremal} means when our $L^*$-base is $CL$ and, as our naming convention suggests, we 
utilize the class of extremally disconnnected frames (see, e.g., \cite[Sec.~III.3.5]{Johnstone:1982} or \cite[Prop.~8.4.1]{BallWW:2002}):

\begin{definition} \label{def: ED}
    A frame $L$ is \emph{extremally disconnected}, abbreviated to ED, if $a^*\join a^{**}=1$ for all $a \in L$, or equivalently, $\B L = CL$.  
\end{definition} 
Note that $L$ is ED iff the non trivial DeMorgan law $(x\meet y)^* =x^* \join y^*$ holds for all $x,y \in L$; equivalently, $x^{**} \join y^{**} = (x^{**} \join y^{**})^{**}$ for all $x,y \in L$. We remind the reader that zero-dimensional frames are completely regular. Taking this into account, results about extremal disconnectivity may be summarized as follows. 

\begin{theorem} \label{thm: ED}\label{thm: left exact}
For $L$ zero-dimensional, the following are equivalent.
\begin{enumerate}[label=\normalfont(\arabic*), ref = \arabic*]
    \item $L$ is $ED$.
    \item $L$ is $CL$-extremal $($that is, $CL$ is a frame$)$.
  \item  $CL= D\!_\infty CL = \B L$.
    \item  $CL$ is a complete Boolean algebra.
 \item The embedding of $CL$ into $L$ has a left adjoint.
    \item $\lambda L$ is $ED$.
     \item $\beta L$ is $ED$.
\end{enumerate}
\end{theorem}
\begin{proof}
It is obvious that $(2)$ and $(4)$ are equivalent. That $(1),(2),(3)$ and $(5)$ are equivalent follows from \cref{thm: A extremal} and \cref{prop: L**}. Note that, since $L$ is zero-dimensional, it follows from \cref{rem: j=k=regular} that the left adjoint in (5), which is given by $k$, is equal to $j$, hence is automatically left exact.\footnote{Alternatively, this can be seen by observing that proving (4) only requires that the embedding of $CL$ into $L$ has a left adjoint.
Indeed, if the embedding $CL \hookrightarrow L$ has a left adjoint $\ell:L\to CL$, then $\bigvee_{CL} S = \ell(\bigvee_L S)$ for any $S \subseteq CL$, and hence $CL$ is a complete Boolean algebra.}
Using these equivalences, (1)$\Leftrightarrow$(6) and (1)$\Leftrightarrow$(7) follow from \cref{CL} that $CL = C(\lambda L) = C(\beta L)$.\footnote{(1)$\Leftrightarrow$(7) was proved using a different method in \cite[Prop.~ 8.4.2]{BallWW:2002}.}
\end{proof}
We now consider $D\!_\sigma CL$, that is, the sub-$\sigma$-frame of $D\!_\infty CL$ generated by $CL$.
 We abuse notation and let $\B A$ denote $\{a^{**}: a\in A\}$ for any subset $A$ of the frame $L$. Unlike $\B L$, $\B A$ need not be a Boolean algebra. Recalling that $L_\sigma$ is the sub-$\sigma$-frame of $L$ generated by $CL$, we have:

\begin{proposition}\label{prop:Dsigma CL}
    For $L$ zero-dimensional, $D\!_\sigma CL= \B L_\sigma$. 
\end{proposition}
\begin{proof}
    By \cref{prop: L**},
  \begin{equation*}
  \begin{split}
        a \in D\!_\sigma CL &\Leftrightarrow a = \displaystyle\bigvee_{D\!_\infty CL} \{c_n: c_n \in CL \} \\
    &\Leftrightarrow a = \displaystyle\bigvee_{\B L} \{c_n: c_n \in CL \} \\
    &\Leftrightarrow a = \left(\displaystyle\bigvee_{L} \{c_n: c_n \in CL \}\right)^{**}\\
    & \Leftrightarrow a \in \B L_\sigma.
      \end{split}
    \end{equation*}
\end{proof}

 To describe the extreme case when $D\!_\sigma CL = L_\sigma$, we utilize a few definitions (some old and some new).
    \begin{definition} \label{def: P-frame} \cite[Def.~8.4.6]{BallWW:2002} 
    For a completely regular frame,
        \begin{enumerate}
            \item $L$ is a  \emph{$P$-frame} if every cozero element is complemented, that is $\Coz L = CL$; 
         \item $L$ is an \emph{almost P-frame} if every cozero element is pseudocomplemented, that is $\Coz L \subseteq \B L $.  \end{enumerate}
    \end{definition}
    \begin{definition} For a zero-dimensional frame,
        \begin{enumerate}
           \item \cite[p.~17]{BB:2014} $L$ is a  \emph{$P_0$-frame} if  the countable join of complemented elements is complemented, that is $ L_{\sigma}=CL$;
           \item $L$ is an \emph{almost $P_0$-frame} if the countable join of complemented elements is pseudocomplemented, that is $L_\sigma \subseteq \B L$.
        \end{enumerate}
    \end{definition}

   \begin{remark} 
    It follows from \cite[Prop. 2]{BB:2014} that (almost) $P$-frames are precisely strongly zero-dimensional (almost) $P_0$-frames. It is an open question as to whether strong zero-dimensionality is required (see \cite[p.~22]{BB:2014}).\footnote{This is related to the open question mentioned after \cref{def: BD}.}
    \end{remark}

\begin{proposition} For $L$ zero-dimensional, $D\!_\sigma CL = L_\sigma$ iff $L$ is an almost $P_0$-frame.
\end{proposition} 
   
\begin{proof} 
Suppose that $D\!_\sigma CL = L_\sigma$. Since $D\!_\sigma CL$ 
is a sub-$\sigma$-frame of 
$D\!_\infty CL = \B L$, we have  
$L_\sigma \subseteq \B L$, that is $L$ is an almost $P_0$-frame. 
Conversely, suppose $L$ is an almost $P_0$-frame, that is $L_\sigma \subseteq \B L$. Then countable joins of elements of $CL$ in $L$ and $\B L$ coincide. But $D\!_\sigma CL$ is countably join-generated in $\B L$ by $C L$, 
while $L_\sigma$ is countably join-generated in $L$ by $CL$. Thus, $D\!_\sigma CL = L_\sigma$.
\end{proof}
Recall from the preliminaries that $L$ is strongly zero-dimensional iff $L_\sigma = \Coz L$, and thus the following is immediate.
    \begin{corollary}
        For  $L$ strongly zero-dimensional, $D\!_\sigma CL = {\sf Coz}\,L$ iff $L$ is an almost P-frame.
  \end{corollary} 

Now we consider when $D\!_\sigma CL = CL$, that is, when $L$ is $CL$-basic. We introduce an obvious weakening of the class of basically disconnected frames.
\begin{definition} \label{def: BD}
A frame $L$ is 
\begin{enumerate}
    \item \cite[Prop.~8.4.3]{BallWW:2002} \emph{basically disconnected}, abbreviated to BD, if $a^* \join a^{**}=1$ for all $a \in \Coz L$, or equivalently, $\B (\Coz L) = CL$;
     and 
    \item  \emph{$BD_0$}  if $a^* \join a^{**}=1$ for all $a \in L_\sigma$, or equivalently, $\B L_\sigma = CL$.
 \end{enumerate}
    \end{definition}
  It is obvious that every BD frame is BD$_0$, and the two notions are equivalent for strongly zero-dimensional frames.
  It remains an open problem whether the two coincide for the larger class of zero-dimensional frames. In fact, this remains an open problem even for spatial $L$. In other words, it is unknown whether each zero-dimensional Hausdorff space $X$ in which the Boolean algebra ${\sf Clop}(X)$ of clopens is $\sigma$-complete must be basically disconnected (see \cite{Dashiell:2013}).
  
    \begin{proposition} \label{prop:BD0} For any zero-dimensional frame,
\begin{enumerate}[label=\normalfont(\arabic*), ref = \arabic*]
    \item $L$ is $BD_0$ iff $CL$ is a sub-$\sigma$-frame of $\B L$.
    \item $L$ is a $P_0$-frame  iff $L$ is almost $P_0$ and $BD_0$.
\end{enumerate}
    \end{proposition}
    
    \begin{proof}
        (1) We have:
        \[
        \begin{array}{llll}
        L \mbox{ is } BD_0 & \Leftrightarrow & \B L_\sigma = CL & \\
        & \Leftrightarrow & D\!_\sigma CL = CL & \mbox{by \cref{prop:Dsigma CL}} \\ & \Leftrightarrow & CL \mbox{ is a sub-$\sigma$-frame of } \B L & \mbox{by \cref{prop: L**}}.
        \end{array}
        \]
        
        (2) If $L$ is a $P_0$-frame, then $L_\sigma = CL$, and thus $\B(L_\sigma) = \B(CL) = CL$, showing $L$ is $BD_0$. But $CL \subseteq \B L$ always, so $L_\sigma \subseteq \B L$, showing that $L$ is almost $P_0$. For the converse, assume $L$ is $BD_0$ and almost $P_0$, and take $a \in L_\sigma$. Then almost $P_0$ gives $a= a^{**}$ and $BD_0$ gives $a^{**} \in CL$, so $L_\sigma = CL$. Hence, $L$ is a $P_0$-frame.
    \end{proof}
We next observe that every completely regular BD frame is strongly zero-dimensional. This can be seen by noting that $L$ is BD iff $\beta L$ is BD, and $\beta L$ BD implies that $\beta L$ is zero-dimensional \cite[Prop.~8.4.4, 8.4.5]{BallWW:2002}. Therefore, $L$ is strongly zero-dimensional (see \cref{prop: Lind corefl3}).
The following is a simple direct proof.

\begin{lemma}
    Every completely regular $BD$ frame is strongly zero-dimensional.
\end{lemma}
\begin{proof} Let $a \in \Coz L$. Then $a =\bigvee\{b_n:b_n \cb b_{n+1}\}$. By \cite[Cor.~1]{BBGilmour:1996},  there are $c_n \in \Coz L$ such that $b_n \cb c_n \cb b_{n+1}$, so $a = \bigvee c_n^{**}$. Basic disconnectivity gives $c_n^{**} \in CL$, so $a \in L_{\sigma}$, showing that $L$ is strongly zero-dimensional.
\end{proof}
The lemma above, together with \cref{prop:BD0} gives the following.
\begin{proposition}\label{prop:BD} 
For $L$ completely regular, 
\begin{enumerate}[label=\normalfont(\arabic*), ref = \arabic*]
    \item $L$ is $BD$ 
 iff $L$ is strongly zero-dimensional and $CL$ is a sub-$\sigma$-frame of $\B L$.
 \item \cite[Prop.~8.4.7]{BallWW:2002} $L$ is a $P$-frame iff $L$ is almost $P$ and $BD$.
   \end{enumerate} 
   \end{proposition}

Utilizing \cref{thm: A basic,prop:Dsigma CL,prop:BD0}, $BD_0$ may be summarized as follows.

\begin{theorem}
    \label{thm: BD_0}
For any zero-dimensional frame $L$, the following are equivalent.
\begin{enumerate}[label=\normalfont(\arabic*), ref = \arabic*]
    \item $L$ is $BD_0$.
    \item $L$ is $CL$-basic $($that is, $CL$ is a $\sigma$-frame$)$.
    \item  $CL = D\!_\sigma CL = \B L_\sigma$.
    \item $CL$ is a $\sigma$-complete Boolean algebra.
\item The embedding of $ CL$ into $ L_\sigma$ has a left adjoint.\footnote{
As in \cref{thm: left exact}, we don't need to require left exact in (5) because proving (4) only requires that the embedding of $CL$ into $L_\sigma$ has a left adjoint. Indeed, if the embedding $CL \hookrightarrow L_\sigma$ has a left adjoint $\ell:L_\sigma\to CL$, then $\bigvee_{CL} a_n = \ell(\bigvee_L a_n)$ for any countable $\{ a_n \} \subseteq CL$. Thus, $CL$ is a $\sigma$-complete Boolean algebra.}
 
\end{enumerate}
\end{theorem}
As a corollary to the above proposition, one gets the following summary of basic disconnectivity for strongly zero-dimensional frames.
\begin{corollary}
    \label{cor: BD}
For any strongly zero-dimensional frame $L$, the following are equivalent.
\begin{enumerate}[label=\normalfont(\arabic*), ref = \arabic*]
    \item $L$ is $BD$.
    \item $L$ is $CL$-basic $($that is, $CL$ is a $\sigma$-frame$)$.
    \item  $CL = D\!_\sigma CL = \B(\Coz L)$.
    \item $CL$ is a $\sigma$-complete Boolean algebra.
    \item The embedding of $ CL$ into $ \Coz L$ has a left adjoint.
    \item $\lambda L$ is $BD$. 
    \item $\beta L$ is $BD$.
\end{enumerate}
\end{corollary}
\begin{proof}
    Recall that for a strongly zero-dimensional frame, $\Coz L = L_\sigma$ and thus the first five equivalences follow from \cref{thm: BD_0}. Since $CL = C(\lambda L) = C (\beta L)$, these five statements are also equivalent to (6) and (7).\footnote{The equivalence (1)$\Leftrightarrow$(7) was proved in \cite[Prop.~8.4.5]{BallWW:2002}.}
\end{proof}

\begin{remark} \label{BDnotED}
There exist BD spaces that are not ED \cite[Example 4N]{GJ:1960}. The frame of opens of any such space is CL-basic but not CL-extremal.
\end{remark}

 We conclude this section by considering the situation when $D\!_\sigma CL = D\!_\infty CL$.
  For strongly zero-dimensional frames, $D\!_\sigma CL = \B(\Coz L)$ (see \cref{prop:Dsigma CL}) and $D\!_\infty CL = \B L$ (see \cref{prop: L**}), so these are frames for which $\B(\Coz L) = \B L$. This is equivalent to the Booleanization map being coz-onto: combining cozero and Boole, such frames were proposed to be called ``coole" in \cite{WW:2020}. (This condition is studied in spaces in \cite{Gruenhage:2006} and \cite{HM:1993} where the terms ``cozero approximated" and ``fraction-dense" are used.) 
  
 \begin{definition}\label{defn:coole} A frame $L$ is \emph{coole} if $\B(\Coz L)= \B L$.
\end{definition}

 \begin{theorem} \label{thm: coole}
 For a strongly zero-dimensional frame $L$, the following are equivalent.
 \begin{enumerate}[label=\normalfont(\arabic*), ref = \arabic*]
 \item $L$ is coole. 
    \item $L$ is $CL$-coole $($that is, $D\!_\sigma CL = D\!_\infty CL)$.
      \item $\B(\Coz L)$ is a frame.
    \item $\lambda L$ is coole. 
    \item $\beta L$ is coole.
    \end{enumerate}    
 \end{theorem}

 \begin{proof} \cref{thm: Dsigma=Dinfty,prop:Dsigma CL,prop: L**} give the equivalence of (1) through (3). Since $L$ and $\lambda L$ always have the same cozeros and Booleanizations (see \cref{fact: known}), it is clear that (1)$\Leftrightarrow$(4). (In fact, this does not require $L$ to be even zero-dimensional.) Now,
     $L$ and $\beta L$ also have the same Booleanization. If $L$ is strongly zero-dimensional, then $\beta L = \I CL$ \cite{BB:1989}. Thus, the two frames have the same complemented elements too, and so the equivalence of (1) and (2) implies the equivalence of (1) and (5). 
     \end{proof}

 To formulate the above for zero-dimensional frames, we introduce the following notion:
 
 \begin{definition}
   A frame $L$ is \emph{0-coole} if $\B L_\sigma = \B L$.  
 \end{definition} 

\cref{thm: coole} has an obvious 0-coole version:
 
 \begin{proposition} For a zero-dimensional frame $L$, the following are equivalent.
 \begin{enumerate}[label=\normalfont(\arabic*), ref = \arabic*]
 \item $L$ is 0-coole.
    \item $L$ is $CL$-coole $($that is, $D\!_\sigma CL = D\!_\infty CL)$.
      \item $\B L_\sigma $ is a frame.
\end{enumerate}    
 \end{proposition}
 
\begin{remark} All the ``0-versions" are equivalent to their ``non-0" counterparts for strongly zero-dimensional frames. 
 As before, it is unclear whether coole and 0-coole coincide for zero-dimensional frames.
 \end{remark}
 
The following  table gives a summary of the situation (the notions of Oz and strongly Oz  are defined in the next section; see \cref{def: Oz}  and \cref{rem:str. Oz}):

\begin{table}[h]\label{table}
\footnotesize
\begin{tabular}{|c|c|c|c|}
\hline
 \textbf{Zero-dimensional} &  $CL$ is base                     &  \textbf{Strongly zero-dimensional} & $CL$ is base                     \\
                                    &                             &                &   and $L_{\sigma} = \Coz L$                        \\ \hline \hline
P$_0$                                            & $L_\sigma = CL$             & P                                                      & $\Coz L = CL$             \\ \hline
almost P$_0$                                     & $L_\sigma \subseteq \B L$ & almost P                                                & $\Coz L \subseteq \B L$ \\ \hline
strongly $Oz$                                    & $\B L\subseteq L_\sigma$  & Oz                                                      & $\B L\subseteq \Coz L$  \\ \hline
BD$_0$                                           & $\B L_\sigma = CL$      & BD                                                      & $\B(\Coz L) = CL$      \\ \hline
                         ED                        &   $\B L = CL$                          & ED                                                      & $\B L = CL$             \\ \hline 
almost P$_0$ and str. Oz      & $\B L=L_\sigma$                             & almost P and Oz (iff ED and P)                                 &  $\B L=\Coz L$                         \\ \hline\hline
P$_0$ iff almost P$_0$ and BD$_0$                &                             & P iff almost P and BD                               &                           \\ 
\hline

\end{tabular}

\medskip
\caption{Summary of zero dimensional vs. strongly zero dimensional.}
    \label{tab:ZD vs str. ZD}
\end{table}

\section{The \texorpdfstring{$L^*$}--base of cozero elements} \label{sec: cozero}

In the preceeding section, we have already seen that cozero elements play a role. We now consider them as our base of interest:  $\Coz L$ is a base iff the frame $L$ is completely regular, and thus for this section we assume all frames to be completely regular. Since $\Coz L$ is a sub-$\sigma$-frame of $L$, it is not only an $L^*$-base, but one such that $\Coz L = D\!_\sigma \Coz L$, and so, for trivial reasons, a frame is completely regular iff it is $\Coz$-basic, and $\Coz$-extremal iff  $\Coz$-coole.  
The two extreme cases for $D\!_\infty \Coz L$ are  $D\!_\infty \Coz L = L$ and $D\!_\infty \Coz L = \B L$.\footnote{
 $D\!_\infty \Coz L$ was studied in \cite{ball_lindelof_2017} where it was denoted by sl($\Coz$).
}
To analyze the interesting situation when \cref{thm: A extremal} holds, that is, when $\Coz L$ is itself a frame (or,  equivalently, $\Coz L = D\!_\infty \Coz L$), we need some definitions.

\begin{definition}\hfill
\begin{enumerate}
    \item A frame is a \emph{cozero frame} if the cozero elements form a frame.
    \item \cite[p.~156]{PicadoPultr:2021} A frame is \emph{perfectly normal} if every element is a cozero, that is $L= \Coz L$.
    \end{enumerate}
\end{definition}
\begin{remark}\label{rem: about CozL}\hfill
 \begin{enumerate}[label=\normalfont(\arabic*), ref = \theremark(\arabic*)]
     \item \label[remark]{rem: about CozL-1}
For any cozero frame $L$, $\Coz L = D\!_\infty \Coz L$, and hence $\Coz L$ is a sublocale of $L$.
By \cref{thm: A extremal},
 the nucleus associated with this sublocale is given by $$x \mapsto \bigwedge \{a \in \Coz L : x \leq a\}.$$
\item Trivially cozero frames are precisely the $\Coz$-coole frames (equivalently, the $\Coz$-extremal frames, as noted above).
\item Every perfectly normal frame is clearly a cozero frame. In general, for a cozero frame $L$, 
joins in $\Coz L$ can be different from those in $L$, and they coincide iff $L$ is perfectly normal. 
(Obviously, meets in $L$ and $\Coz L$ coincide since $\Coz L$ is a sublocale of $L$.) For a cozero frame that is not perfectly normal, see \cref{rem: almost Boolean 2} or \cref{example: Q}. 
\end{enumerate}
\end{remark}

 By \cref{thm: A extremal}, we get the following extreme cases: 

 \begin{proposition}\ \label{prop: cozero frames} 
     \begin{enumerate}[label=\normalfont(\arabic*), ref = \theproposition(\arabic*)]
         \item $\Coz L =D\!_\infty \Coz L = L$ iff $L$ is perfectly normal; 
         \item \label[proposition]{cozero frames-2} $\Coz L = D\!_\infty \Coz L = \B L$ iff $L$ is an extremally disconnected P-frame.
     \end{enumerate}
 \end{proposition}
 \begin{proof}
    The first result is straightforward. For (2), the reverse implication follows from the definitions, and the other follows from the observation that $\Coz L = \B L$ forces  that the binary join in $\B L$ is the same as that in $L$, and hence {$CL = \B L (=\Coz L)$}.
 \end{proof}

 \begin{remark}\ \label{rem: almost Boolean}
 \begin{enumerate}[label=\normalfont(\arabic*), ref = \theremark(\arabic*)]
     \item \label[remark]{rem: almost Boolean 1} It is often noted, without proof, that extremally disconnected P-frames are precisely those for which $\Coz L  = \B L$; we give the argument above for the sake of completeness.

\item \label[remark]{rem: almost Boolean 2} Banaschewski \cite[Sec.~2]{banaschewski_nonmeasurable_2015} 
defines \emph{almost Boolean} frames and shows that these are precisely the extremally disconnected P-frames, thus any non-Boolean, almost Boolean frame is a cozero frame that is not perfectly normal. It is straightforward to see that these are also the almost P, Oz frames (see \cref{tab:ZD vs str. ZD}). The latter are defined below (see \cref{def: Oz}).
 \end{enumerate}
\end{remark}

Recall that $\lambda L$, the Lindel\"of coreflection of a completely regular frame $L$, has the same cozero part as $L$ (see \cref{CozL}), and so the next lemma is obvious.
\begin{lemma} $L$ is a cozero frame iff $\lambda L$ is a cozero frame.
\end{lemma}

The same is not true for perfect normality: in fact for Lindel\"of frames, perfect normality is equivalent to hereditary Lindel\"of, and so 
the Lindel\"of coreflection is perfectly normal iff the frame itself is Lindel\"of \cite{WW:2020}. Thus, any proper Lindel\"of coreflection is never perfectly normal, and we obtain:

\begin{proposition} \label{prop: cozero not pn}
The Lindel\"of coreflection of a non-Lindel\"of cozero frame is cozero but not perfectly normal.
\end{proposition} 

The above proposition is a source of interesting examples.
\begin{example} \label{example: Q}
    \hfill
    \begin{enumerate}[label=\normalfont(\arabic*), ref = \theexample(\arabic*)]
        \item \label[example]{ex: Q1} Since a metric space is separable iff it is Lindel\"of \cite[Cor.~4.1.16]{Eng:1989}, the Lindel\"of coreflection of the frame of opens of any non-separable metric space is cozero but not perfectly normal by \cref{prop: cozero not pn}.
      If, in addition,  $X$ is zero-dimensional, but not discrete
      (e.g., $\mathbb R^\omega$ with the Baire metric \cite[Example 4.2.12]{Eng:1989}), 
      then the clopens of $X$ are properly contained in the cozeros of $X$. Therefore, by \cref{prop: cozero not pn} and \cref{fact: known}, the Lindel\"of coreflection of the frame of opens of $X$ is a cozero frame $L$ such that 
        $CL \subsetneqq \Coz L \subsetneqq L$. Note that this cozero frame is not spatial. 
        \item For a spatial example, take a perfectly normal space $X$ that is not realcompact. That such an $X$ exists follows from \cite{Ost:1976a}, although the proof is not within ZFC. Let $L$ be the frame of opens of the realcompactification of $X$. Then $L$ is a spatial cozero frame that is not perfectly normal. (For details about realcompactness see, e.g., \cite[Sec.~8]{GJ:1960}.)
        \end{enumerate}
\end{example}

To further describe cozero frames, we consider where this class of frames sits relative to other classes. We already know that it properly includes all perfectly normal frames. 
\begin{definition}\label{def: Oz}
 \cite[Prop.~2.2]{BDGWW:2009} A frame is called \emph{Oz} if every pseudocomplement is cozero, that is $\B L\subseteq \Coz L$.
 \end{definition}

\begin{remark}\label{rem:str. Oz}
    A natural idea is the ``0" version of Oz: A frame is called \emph{strongly Oz} if every pseudocomplement is a countable join of complemented elements (that is, $\B L\subseteq L_\sigma$).
Strongly Oz clearly
 implies Oz,  and thus it is a more appropriate name than Oz$_0$ (all the other ``0-notions" studied in \cref{sec: zero-dimensional}  are weaker than their Coz-counterparts).
One can consider various interplays between these notions. For example, we saw above that $\Coz L = \B L$ is equivalent to almost P and Oz (iff P and ED iff almost Boolean). Similarly, $L\!_\sigma = \B L$
     is equivalent to almost P$_0$ and strongly Oz  but it is unclear what this means for zero-dimensional frames that are not strongly zero dimensional. 
     
\end{remark}

Obviously perfect normality implies Oz, and so does being cozero:

\begin{proposition} Every cozero frame is Oz.
\end{proposition}
\begin{proof} 
Recall that $D\!_\infty \Coz L$ is a dense sublocale of $L$ and so, as sublocales, $\B L \subseteq D\!_\infty \Coz L$. But for cozero frames the latter is $\Coz L$ (see \cref{rem: about CozL-1}).
\end{proof}
ED frames are a proper subclass of Oz frames (for example, the frame of opens of $\mathbb R$ is Oz but not ED), but not all ED frames are cozero as the following results clarify.

\begin{proposition} \label{prop: finite B}For any $\sigma$-complete Boolean algebra $B$, the following are equivalent.
\begin{enumerate}[label=\normalfont(\arabic*), ref = \arabic*]
\item Every ideal of $B$ is countably generated.
\item The countably generated ideals of $B$ form a frame.
\item $B$ is finite.
\end{enumerate}
    \end{proposition}
    
  \begin{proof}
    This is a consequence of the following proposition, which is of independent interest.
   \end{proof}

\begin{proposition}\label{prop: infinite-Boolean} Every infinite $\sigma$-complete Boolean algebra $B$ has an ideal which is not countably generated.
\end{proposition}
    
    \begin{proof} 
 Let $X$ be a countably infinite family of pairwise disjoint elements of $B$ (see, e.g., \cite[Prop.~3.4]{Koppelberg:1989}). Let $\{X_\alpha :\alpha\in \Gamma\}$ be an uncountable family of pairwise almost disjoint subsets of $X$ (see, e.g., \cite[6Q.1]{GJ:1960}). Since $B$ is $\sigma$-complete, each $X_\alpha$ has a join $b_\alpha$ in $B$. Let $I$ be the ideal of $B$ generated by $\{b_\alpha : \alpha \in \Gamma \}$.\   We show that $I$ is not countably generated. 
 
 First, we claim that for any $\alpha\in\Gamma$, $b_\alpha$ is not in the ideal generated by $\{b_\sigma :\sigma\neq\alpha\}$. Otherwise, there exist $\sigma_1\ldots\sigma_n$ with $b_\alpha\leq \bigvee_{i=1}^n b_{\sigma_i}$. But this is impossible because $X_\alpha$ is almost disjoint from the $X_{\sigma_i}$ so there exists $x\in X_\alpha\setminus\bigcup_{i=1}^n X_{\sigma_i}$. Since $x\wedge  y=0$ for all $y\in X_{\sigma_i}$, $x\wedge b_{\sigma_i}=0$ by 
 $\sigma$-distributivity in the Boolean algebra $B$. Because $x\le b_\alpha$,  it is not possible that $b_\alpha\leq \bigvee_{i=1}^n b_{\sigma_i}$. 
 Thus, no countable subset of the $b_\alpha$ can generate $I$. 
 
 Lastly, for an arbitrary countable family $\{a_n\}\subseteq I$, each $a_n$ is
 underneath a finite join from $\{b_\alpha:\alpha\in\Gamma\}$. If $\{a_n\}$ generated $I$, this would provide a countable generating subset of $\{b_\alpha\}$, which is a contradiction.
    \end{proof}

\begin{remark} 
Since $\sigma$-complete Boolean algebras correspond to basically disconnected compact Hausdorff spaces (see, e.g., \cite[Prop.~7.21]{Koppelberg:1989}), the preceding proposition yields: an infinite basically disconnected compact Hausdorff space (or compact regular frame) is never perfectly normal.
\end{remark}
    
    \begin{example} \label{ex: betaN}
 The frame $L$ of opens of $\beta\mathbb{N}$ is isomorphic to $\I \wp\mathbb{N}$, the ideal frame of the powerset of $\mathbb{N}$. Therefore, $L$ is an extremally disconnected frame. It is not a cozero frame by \cref{prop: finite B} because $\Coz \I \wp\mathbb{N}$ is precisely the countably generated ideals of $ \wp\mathbb{N}$. This yields an example of an ED frame that is not cozero, and a spatial one at that.
   \end{example} 
  In \cref{BDnotED} we noted that ED is stronger than BD. However, if $L$ is a cozero frame the difference collapses.

\begin{lemma}    
A $BD$ cozero frame is $ED$. 
\end{lemma} 

\begin{proof} 
$L$ being BD gives $\B(\Coz L)=CL$. Since $L$ is also a cozero frame, $\Coz L$ is a dense sublocale of $L$, and thus, as sublocales, $\B L \subseteq \Coz L$. Therefore, applying the nucleus $(-)^{**}$ gives  $\B L= \B (\Coz L) = CL$. 
\end{proof}

\cref{ex: betaN} shows that the converse is not true.

\begin{remark}\

\begin{enumerate}
    \item There are some obvious connections between cozero frames and the notions of tightness (see \cite{ball_tightness_2016}), which will be explored in future work.\footnote{A note of caution: in \cite{ball_tightness_2016} 
 cozero means perfectly normal.}
 \item  \label{rem: non-D-bases}
Another source of interesting bases are the ``cozero bases''  for completely regular frames, and their zero-dimensional analogs ``c-bases'', both defined in \cite{BBGilmour:2001}. The latter are frame theoretic analogs of Boolean bases in topology \cite{Dwinger:1961}. These are always L$^*$-bases, as are the uniformly cozero elements of a uniform frame (see, e.g., \cite{Walters:1989}). One can clearly obtain corresponding notions to those laid out in \cref{sec: zero-dimensional,sec: cozero} for these bases.
\end{enumerate}
\end{remark}

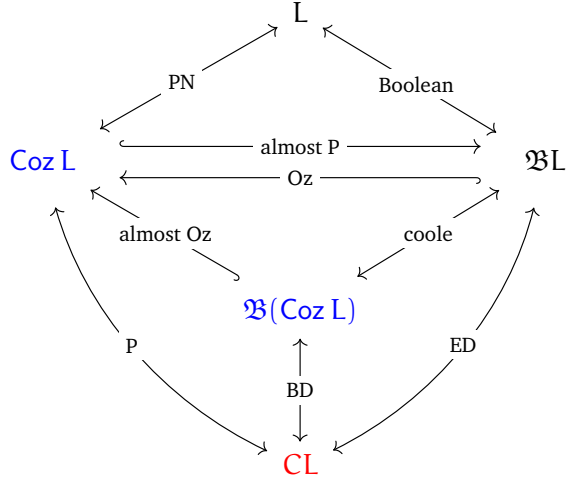
\begin{figure}[h] 
\[\begin{tikzcd}
	&& L \\
	\\
	\color{blue}{\Coz L} &&&& \B L \\
	\\
	&& \color{blue}{\B (\Coz L)} \\
	\\
	&& \color{red}{CL}
    \arrow["{\scriptsize\text{PN}}"{description}, shorten <=6pt, tail reversed, from=3-1, to=1-3]
	\arrow["{\scriptsize\text{Boolean}}"{description}, shorten >=6pt, tail reversed, from=1-3, to=3-5]
        	\arrow["{\scriptsize\text{P}}"{description}, curve={height=24pt}, shorten <=10pt, tail reversed, from=3-1, to=7-3]
	\arrow["{\scriptsize\text{ED}}"{description}, curve={height=-24pt}, shorten <=10pt, tail reversed, from=3-5, to=7-3]
	\arrow["{\scriptsize\text{BD}}"{description}, tail reversed, from=5-3, to=7-3]
	\arrow["{\scriptsize\text{coole}}"{description}, shorten <=6pt, shorten >=6pt, tail reversed, from=5-3, to=3-5]
	\arrow["{\scriptsize\text{almost P}}"{description}, shift left=2, shorten <=11pt, shorten >=11pt, hook, from=3-1, to=3-5]
	\arrow["{\scriptsize\text{Oz}}"{description}, shift left=2, shorten <=11pt, shorten >=11pt, hook, from=3-5, to=3-1]
	\arrow["{\scriptsize\text{almost Oz}}"{description}, shorten <=6pt, shorten >=6pt, hook, from=5-3, to=3-1]
\end{tikzcd}\]
\caption{Summary for completely regular frames.}
    \label{fig: Coz summary}
\end{figure}

\cref{fig: Coz summary} summarizes the situation for completely regular frames.
  We recall that a frame is \emph{weak Oz} if the pseudocomplement of a cozero is cozero (see \cite[Sec.~5]{BDGWW:2009}) and a slightly weaker notion is that of \emph{almost Oz}: the double pseudocomplement of a cozero is cozero, that is $\B (\Coz L) \subseteq \Coz L$ (see \cite{WW:2020}). We include this here for the sake of ``symmetry" of the diagram. 

There is obviously a similar diagram with $L_\sigma$ replacing $\Coz L$ and all the analogous ``0-versions'' (see \cref{tab:ZD vs str. ZD}) with some potentially new classes to be considered.

\section{The \texorpdfstring{$L^*$}--base of compact elements} \label{sec: compact}

For any frame $L$, consider the poset of compact elements of $L$, denoted by $KL$.

\begin{definition} \cite[pp.~63--64]{Johnstone:1982} 
$L$ is a \emph{coherent frame} if $KL$ is a bounded sublattice and a base of $L$.     
\end{definition}

Coherent frames are also examples of frames with an $L^*$-base.
As documented in Johnstone \cite[p.~64]{Johnstone:1982}, $L$ is a coherent frame iff $L$ is isomorphic to the frame of ideals of a bounded distributive lattice, namely $L \cong \I KL$.
 Since $KL$ is bounded, the two extreme cases for $D\!_\infty KL$ are $D\!_\infty KL = \B L$ and $D\!_\infty KL = L$.
\begin{proposition} \label{thm: coherent}
    Let $L$ be a coherent frame.
    \begin{enumerate}[label=\normalfont(\arabic*), ref = \arabic*]
        \item $D\!_\infty KL = \B L$ iff $KL \subseteq \B L$.
        \item $D\!_\infty KL = L$ iff every ideal of $KL$ is a D-ideal.
    \end{enumerate}
\end{proposition}
\begin{proof} 
(1) The left-to-right implication is clear. For the other implication, we show that $KL$ is a base for $\B L$: Let $a\in \B L$. Since $KL$ is a base for $L$, $a=\bigvee S$ for some $S\subseteq KL$. Therefore, $a=\left(\bigvee S\right)^{**}$, so $a=\bigvee_{\B L} S$ (because $S\subseteq KL\subseteq\B L$). Thus,  
$KL$ is an S-base for $\B L$, and so $D\!_\infty KL \subseteq \B L$ by \cref{thm: interval}. The other inclusion follows since $KL$ bounded implies $D\!_\infty KL$ is a dense sublocale of $L$, and hence contains $\B L$. 

(2) Since $KL$ is a bounded distributive lattice, every D-ideal of $KL$ is an ideal. Therefore, every ideal of $KL$ being a D-ideal is equivalent to $\I KL = D\!_\infty KL$. Thus, the result follows since $L \cong \I KL$.
    \end{proof}

To analyze $D\!_\sigma KL$ we work with coherent $\sigma$-frames, the $\sigma$-frame analogs of coherent frames (see \cite{Walters:1989}). Observe that $L_{\sigma_{KL}}$ is the $\sigma$-frame of Lindel\"of elements of $L$, which we denote by $K\!_\sigma L$. 
Since $L\cong \I KL$, we obtain that 
$K\!_\sigma L$
is isomorphic to the coherent $\sigma$-frame of countably generated ideals of $KL$, which we denote by $\I_\sigma KL$. The two extreme cases become $D_\sigma KL = \B(K\!_\sigma L)$ and $D_\sigma KL = \I_\sigma KL$.

\begin{proposition} \label{thm:}
    Let $L$ be a coherent frame.
    \begin{enumerate}[label=\normalfont(\arabic*), ref = \arabic*]
        \item $D\!_\sigma KL = \B(K\!_\sigma L)$ iff $KL \subseteq \B(K\!_\sigma L)$.
        \item $ D\!_\sigma KL = \I_\sigma KL$ iff every countably generated ideal of $KL$ is a countably generated D-ideal. 
    \end{enumerate}
\end{proposition}

\begin{proof} 
The proof is an appropriate adjustment of that of \cref{thm: coherent} by replacing arbitrary joins with countable joins and ideals by countably generated ideals. 
\end{proof}
Applying \cref{thm: A extremal}, \cref{thm: A basic}, and \cref{thm: Dsigma=Dinfty}, we  obtain:

\begin{theorem} \label{thm: KL-extremal}
    Let $L$ be a coherent frame.
    \begin{enumerate}[label=\normalfont(\arabic*), ref = \arabic*]
        \item $L$ is $KL$-extremal iff every D-ideal of $KL$ is principal.
        \item  $L$ is $KL$-basic iff every countably generated D-ideal of $KL$ is principal.
        \item $L$ is $KL$-coole iff every D-ideal of $KL$ is countably generated.
    \end{enumerate}
\end{theorem}

Of course one can get all the other equivalent conditions listed in these theorems. 
\begin{remark}\hfill
\begin{enumerate}
    \item $KL$-extremal frames were first considered in \cite[Sec.~9]{BH14} under the name of extremally disconnected stably compact frames. 
    \item It is informative to compare the coherent case to the zero-dimensional case. 
    In particular, in the future we aim to develop the subclasses of coherent frames similar to the ones we developed for the zero-dimensional frames (see \cref{tab:ZD vs str. ZD}  and \cref{fig: Coz summary}). 
   \end{enumerate}
\end{remark} 

\section{Conclusion}

As we have already highlighted, there are plenty of unresolved issues even for bases that are fairly well understood, and some interesting new ideas have come to light, notably the non-trivial notion of a cozero frame. Moreover, applying the machinery developed in \cref{sec: S-bases,sec: completeness,sec: D-bases} to other classes of frames with specific bases promises to be fruitful. For example, 
 we recall that a frame $L$ is {\em algebraic} if $KL$ is a base for $L$. An algebraic frame is {\em arithmetic} (\cite{Compendium2003}) or an {\em M-frame} (\cite{IberkleidMcGovern2009,Bhattacharjee2019}) if $a,b\in KL \Rightarrow a\wedge b\in KL$.  Clearly, $L$ is arithmetic iff $KL$ is an L-base and coherent iff $KL$ is an L$^*$-base. 

 Much of what we discussed in \cref{sec: compact} generalizes to the setting of arithmetic frames. 
However, our results do not generalize directly to all algebraic frames since $KL$ is not an S-base of an algebraic frame $L$. We leave it to future work to develop appropriate machinery to handle all algebraic frames.

It would be interesting to generalize $D\!_\sigma$ to
 $D\!_\kappa$ (\`a la Madden \cite{Madden:1995}), where $\kappa$ is an arbitrary infinite cardinal, thus providing $\kappa$ variations of BD and BD$_0$.
 Note that BD and BD$_0$ are really ``countable'' weakenings of ED. Thus, one would obtain the notions of A-$\kappa$-extremal and DA-$\kappa$-extremal, where $A$-basic is A-$\sigma$-extremal etc. It seems plausible that this would give subtle ways to better understand the complicated structure of frames (and spaces). 

 Another interesting direction is to consider a version of \cref{thm: interval} for $\sigma$-frames or, more generally, ``topless'' $\sigma$-frames. The theory of the latter still needs to be developed. 

We suggest the following (incomplete) list of special elements that (promise to form interesting bases and) should be explored and viewed through the lens we have developed in this paper. We give a smattering of results that are already known and leave the rest to future work.
\begin{enumerate}
    \item Consider $\B(\Coz L)$: if $L$ is completely regular then $\B(\Coz L)$ is an S$^*$-base, and it is an L$^*$-base  iff   $x^{**} \join y^{**} = (x^{**} \join y^{**})^{**}$ for all $x,y \in \Coz L$; i.e., $L$ is BD.
    \item The collection $\mathrm{Pel}(L) = \{ a\in L : x\in\Coz L \mbox{ and } x \le a \Rightarrow x \prec a\}$ of all P-elements of a completely regular frame $L$ studied in \cite{Dube:2019} is an S$^*$-base of $L$. 
If $L$ is BD, then $\mathrm{Pel}(L)$ is a dense sublocale of $L$ \cite[Thm.~4.1]{Dube:2019}. Moreover, the two extreme cases are: $\mathrm{Pel}(L)=L$ iff $L$ is a P-frame (straightforward) and $\mathrm{Pel}(L)=\B  L$ iff the P-sublocales of $L$ form a frame \cite[Thm.~4.3]{Dube:2019}.
 \item For an arithmetic frame $L$, the collection $dL = \{ a\in L : a =\bigvee\{k^{**} : \ k\in KL, k \le x  \} \}$ of d-elements is a dense sublocale of $L$ such that $\B (KL)$ is an S-base of $dL$  \cite[p.~247]{MartinezZenk2003}. 
\end{enumerate}

\bibliographystyle{amsplain}

\bibliography{references.bib}

\end{document}